 \magnification=1200
 \pretolerance=999 \tolerance=1999
\overfullrule=0pt    \parskip=20pt \parindent=25pt

\vsize=560pt
\font\blackb=msbm10
\font\gothic=eufm10

\font\bfg=cmbx12

\def\pph#1{p_{#1}^{\vphantom{(1)}}}       
\def\dfrac#1#2{{\displaystyle{#1\over#2}}}
\def\gothi#1{\hbox{\gothic{#1}}}
\overfullrule=0pt

\vglue 25pt

\centerline{\bfg Painlev\'e-type differential equations for the recurrence}
\centerline{\bfg coefficients of semi-classical orthogonal polynomials.}

\bigskip
\centerline{Alphonse P. Magnus}

\centerline{ Institut Math\'ematique,  Universit\'e Catholique de Louvain}

\centerline{Chemin du Cyclotron 2}

\centerline { B-1348 Louvain-la-Neuve}

\centerline{ Belgium}

\centerline{E-mail: {\tt magnus@anma.ucl.ac.be}}

\bigskip

\vbox{
\leftskip=60pt \rightskip=30pt
\noindent{\bf Abstract.}

\noindent
Recurrence coefficients of semi-classical
 orthogonal polynomials (orthogonal
 polynomials related to a weight function $w$ such that $w'/w$ is a
 rational function) are shown to be solutions of non linear differential
 equations with respect to a well-chosen parameter, according to principles
 established by D.\& G.~Chudnovsky.
  Examples are given.
For instance, the recurrence coefficients in $a_{n+1}p_{n+1}(x)=xp_n(x)
-a_np_{n-1}(x)$ of the orthogonal polynomials related to the weight
$\exp(-x^4/4-tx^2)$ on {\blackb R\/} satisfy
$4a_n^3\ddot a_n = (3a_n^4+2ta_n^2-n)(a_n^4+2ta_n^2+n)$, and $a_n^2$
satisfies a Painlev\'e ${\rm P}_{\rm IV}$ equation.
     }

{\bf 1. Introduction: measures and recurrence coefficients of
       orthogonal polynomials.}

Let $\{p_n\}_0^\infty$ be the set of orthonormal polynomials related to
some measure $d\mu$ on its support $S$:
   $$\int_S p_n(x)p_m(x)\,d\mu(x)=\delta_{m,n}.  \eqno(1)$$
The most remarkable property of the $p_n$'s is
the recurrence relation joining them:
 $$ a_{n+1}p_{n+1}(x)= (x-b_n)p_n(x)-a_np_{n-1}(x).  \eqno(2) $$
An often encountered problem in applied and numerical mathematics as well
as in physics is then to relate the coefficients $a_n$ and $b_n$ of
(2) to properties of the measure $d\mu$.

For instance, interesting solid-state
Hamiltonian operators submitted to the
so-called ``recursion method'' (or Lanczos method) show a
tri-diagonal matrix representation. Investigation of spectral properties
of the operator is then equivalent to investigating the measure of (1)
from the recurrence coefficients of (2) [GaCL] [Hay] [HayN]
 [LaG] [LiMu] [OW].

The study of
special partition functions in statistical physics and
quantum physics leads to relations which can be translated as
properties of particular recurrence coefficients. Much important
work is currently done on this subject [Bes] [BIZ] [Fok1] [Fok2]
[Fra1] [Fra2] [GrM1] [GrM2] [HH] [KvM] [LW] [Mo] [Y] [Zu].

Numerical implementation of spectral methods and quadrature formulas
needs accurate determination
of recurrence coefficients for various measures. This appears in the
survey  [Gau] and in some recent papers as [Chin] and [ClS] (see also
the references in [BeR]).

Quite a number of theoretical studies have appeared on this problem of
relating properties of the orthogonality measure to the recurrence
coefficients, especially to their asymptotic behaviour. See at least
the books [Chi], [Fr1], [VA] and the survey papers [Lub] and [GFOPCF].

To give just a taste of the matter, the asymptotic behaviour of the
recurrence coefficients associated to $d\mu(x)=|x|^\rho\exp(-x^4)dx$ on
$x\in${\blackb R} appears in an amazing number of fields:

\item{1.} This extension of Hermite polynomials is studied by
    Shohat [Sho], using methods of Laguerre [Lag]. Later on, Freud [Fr2]
    rediscovered Shohat's formulas (see (31) with $t=0$) and proved that
    $a_n\sim (n/12)^{1/4}$ when $n\to\infty$. Much more has been done
    since then [Lub] [Mag2] [Mag3] [Nev] [GFOPCF] [Nev2], it has been
    shown that a behaviour $d\mu(x)\sim\exp(-|x|^\alpha)dx$ for
    $x\to\pm\infty$ implies $a_n$ and $b_n\sim$ constant $n^{1/\alpha}$ for
    large $n$.

\item{2.} Similar exponential weights were encountered in solid-state and
    statistical physics, where the same asymptotic connections have been
    used (sometimes after clever guesswork) [LiMu] [OW].

\item{3.} These extensions of Hermite and Laguerre polynomials
       also appear in numerical quadrature methods intended
        to solve   Boltzmann and
       Fokker-Planck equations [ClS], where they are called ``speed'',
       ``bimode'' and ``Druyvesteyn'' polynomials.  The recurrence
      coefficients can be computed safely from a suitable algorithmic
      use of the Shohat-Freud
      equations [LeQ], or from asymptotic expansions [ClS]
      (see [Nev2] pp.462-463).

\item{4.} The same orthogonal polynomials reappear in special solutions
     of important differential equations of mathematical physics
      ([Bes] [KvM]; Shohat-Freud's equations are called ``discrete
       Painlev\'e equations'' in [Fok2]), as
     well as in continued fraction expansions of special irrational
     numbers [Chu2].

Any advance in one of these fields is liable to benefit to the
other ones, although the dialog is not always obvious: established
theorems may sometimes have poor constructive contents and be unable to
inspire valuable algorithms; explicit formulas  (using for instance
exotic special functions or high-order determinants) may be delightful
solutions for some people and useless nightmares for other ones;
successful numerical tricks or self-consistent ``{\it Ansatze\/}''
may be out of reach of contemporary methods of proof, etc.

\smallskip

Let us return now to the problem considered here: to deduce properties
of the recurrence coefficients $a_n$ and $b_n$ from the measure
$d\mu(x)$. The Chebyshev orthogonal polynomials are related to measures
involving the square root of a polynomial of degree 2 and yield
constant recurrence coefficients (the simplest case!).
 The classical orthogonal polynomials (Jacobi polynomials and
their limit cases) have a known measure and known recurrence coefficients
($a_n^2$ and $b_n$ are special rational functions of $n$). We may define
a more general class by allowing $a_n^2$ and $b_n$ to be general
rational functions of $n$ (Pollaczek class [Chi]) but then the orthogonality
measure becomes difficult to control. Natural extensions of Chebyshev
polynomials are related to measures involving the square root of a
polynomial of degree $>2$. One finds then an oscillatory behaviour
of the recurrence coefficients
([GV1] [GV2] [Gr] [I] [VA]), whose description may even
need elliptic or hyperelliptic functions [Ak] [Apt] [GaN]. We will
deal here with a further extension,
the so-called semi-classical class (to be defined in the
next section) which represents still a reasonable trade-off between
measure description (easy and containing useful cases) and the
possibility of
description of recurrence coefficients (interesting nonlinear
relations). Much of the work was already done in the end of the
nineteenth century by Laguerre [Lag] who recognized (in 1885) that special
cases (worked by Jacobi and Borchardt)
would indeed involve elliptic functions. As he could not establish the
general recurrence coefficients behaviour, we could suspect that
{\it special functions still unknown in the nineteenth century\/}
would be needed\dots  Painlev\'e transcendents will indeed appear,
and they were investigated in the early twentieth century (see the foreword
of [Pain]).

\medskip

{\bf 2. Formal semi-classical orthogonal polynomials.}

\smallskip

Orthogonal polynomials $p_n$ are usually defined through a measure, so
to satisfy (1).
The construction of these polynomials only needs the sequence
of moments $\mu_k = \int_S x^k d\mu(x), k=0,1,\ldots$
{\it Formal\/} orthogonal
polynomials are only related to a numerical (real or complex) sequence
of numbers $\mu_k, k=0,1,\ldots$, ignoring whether these numbers are
actual moments of some weight or distribution on some support or not.
The polynomial $p_n(z) = \gamma_n z^n +\gamma_{n,1}z^{n-1}+\cdots
+p_n(0)$ is then obtained from the equations
$\gamma_n\mu_{n+k}+\gamma_{n,1}\mu_{n+k-1}+\cdots+p_n(0)\mu_k=0,
k=0,1,\ldots,n-1$ and
$\gamma_n[\gamma_n\mu_{2n}+\gamma_{n,1}\mu_{2n-1}+\cdots+
     p_n(0)\mu_n]=1$. These equations can be solved for any $n=0,1,\ldots$
if the Hankel determinants built with $\mu_0,\ldots,\mu_{2n}$ do not
vanish ([Bre], [dBvR] \S\ 7 ,     see the definition of  {\it regular\/}
formal orthogonal polynomials on p.47 of [Dra] \S\ 1.1-1.3).

If we define a linear form ${\cal L}$ on the space of polynomials by
${\cal L}(x^n)=\mu_n, n=0,1,\ldots$, the polynomials $p_n$
satisfy ${\cal L}(p_np_m)=\delta_{n,m}, n,m=0,1,\ldots$ ([Mar], [Mar2], where
${\cal L}$ is written ${\cal L}_0$).

Regular formal orthogonal polynomials always satisfy the
{\it recurrence relation\/} (2), with
$p_0=\gamma_0=1/\sqrt{\mu_0}, \qquad a_1 p_1(z) = (z-b_0)p_0(z)$,
and where
 $b_0=-\gamma_{1,1}/\gamma_1$, $a_n=\gamma_{n-1}/\gamma_n$,
  $b_n = \gamma_{n,1}/\gamma_n - \gamma_{n+1,1}/\gamma_{n+1},
   n=1,2,\ldots$ ([Bre]  ,  [Dra] \S\ 1.4).

By introducing the formal series
$$f(z)=\sum_0^\infty \mu_k z^{-k-1},  \eqno(3)$$
the equations for $p_n$ are summarized as
 $$ f(z)\pph n(z) = p_{n-1}^{(1)}(z) + \varepsilon_n(z),
    \qquad \varepsilon_n(z) = \gamma_n^{-1}z^{-n-1}+O(z^{-n-2}),
  \eqno(4)  $$
where $p_{n-1}^{(1)}$ is a polynomial of degree $n-1$, (associated
polynomial to $p_n$). These polynomials, as well as the $\varepsilon_n$'s,
satisfy the same recurrence relations (2), but with $p_{-1}^{(1)}=0$,
$p_0^{(1)}=\mu_0\gamma_1=1/(a_1\gamma_0)$. The following relation
$$ \pph n p_{n-2}^{(1)} - \pph{n-1} p_{n-1}^{(1)} =
  \pph{n-1}\varepsilon_n -\pph n \varepsilon_{n-1} =
                                         -1/a_n \eqno(5)$$
is well known ([Chi] , [Fr1] , etc.) From the recurrence relations (2), we have
the main terms in the expansions of $p_n$ and $\varepsilon_n$, which
will be useful later:
$$ p_n(z)= \gamma_n\left[z^n - \left( \sum_0^{n-1} b_i \right) z^{n-1}
         +\left( \sum_{i<j<n}b_ib_j -\sum_1^{n-1}a_i^2\right)z^{n-2}
          +\cdots \right] \eqno(6)$$
$$ \varepsilon_n(z)=
 \gamma_n^{-1}\left[z^{-n-1} + \left( \sum_0^{n} b_i \right) z^{-n-2}
       +\left( \sum_{i\le j\le n}b_ib_j +\sum_1^{n+1}a_i^2\right)z^{-n-3}
          +\cdots \right] \eqno(7)$$
(for the latter one, use $\gamma_n\varepsilon_n(z) =
   (z-b_n)^{-1}\gamma_{n-1}\varepsilon_{n-1}(z) +
   (z-b_n)^{-1}a_{n+1}^2\gamma_{n+1}\varepsilon_{n+1}(z)$).

Of course, if we happen to know a true function of the complex variable
$z$ having the asymptotic expansion (3) when $z\to\infty$ in some way,
and if this function is analytic outside a set $S$ made of contours and
arcs, we may use a Cauchy-like integral representation
$$ f(z) = \int_S w(x) (z-x)^{-1} dx\,, \qquad z\notin S \eqno(8) $$
allowing to recover the convenient description in terms of a
``weight function'' $w$, but the description is not unique and $w$ may
be complex. We then have an integral  representation
of the form ${\cal L}$: ${\cal L}\varphi=\int_S \varphi(x)w(x)dx$.
For instance the Bessel orthogonal polynomials are defined
by $\mu_n = 1/n!, n=0,1,\ldots$ and can be considered as orthogonal
with respect to the complex weight $(2\pi i)^{-1}\exp x^{-1}$ on any
contour containing the origin in its interior.    Remark that orthogonality
of two complex function $\varphi$ and $\psi$
always involves here the product $\varphi\psi$
and {\it not\/} the product $\varphi\overline\psi$ (as in [StT]).

For an example showing how formal orthogonal polynomials can be
investigated through their generating function of formal moments (3),
consider $f(z) = [A(z)-B(z)^{1/2}]/C(z)$, where $A,B$ and $C$ are given
polynomials ({\it formes du second degr\'e\/} in [Mar2] p.122, Def.~7.4).
 Such a function can be represented as (8)
 outside a systems of cuts $S$ joining the zeros of $B$ in some way.
 Here $w(x)$ will have the form
 $w(x)=
\pm(\pi i)^{-1}B(x)^{1/2}/C(x)$ ([N]\S 1.2 \& 4.3.1). If $B$ has only real
zeros, this is a way to introduce special orthogonal polynomials
on several intervals (the intervals where $B(x)\le 0$). Now, (4) gives
here
$ -B^{1/2}p_n = q_n +C\varepsilon_n$,
with $q_n= -Ap_n +Cp_{n-1}^{(1)}$. Squaring yields
$Bp_n^2-q_n^2 = L_n$, where $L_n$ {\it must\/} be a polynomial of
degree bounded by a constant, as the left-hand side is a polynomial, and
as the right-hand side is
 $2q_nC\varepsilon_n+C^2\varepsilon_n^2$. So, $p_n$
is such that the square of this polynomial times a given polynomial $B$
equals the square of another polynomial plus a polynomial of bounded
degree. This is enough for experts to describe $p_n$ in terms of
(hyper)elliptic function and integrals, theta functions, etc. (see
[Ak] \S 53, [Apt], [Brez] pp.\ 296-298,
[N]\S 4.3 , [Peh] ), and to discuss periodic features
in the sequence of the recurrence coefficients ([GV1], [GV2], [Gr], [I],
[Peh1] ). For arithmetic continued fractions connected to Pell's
equation, see [Brez] pp.\ 39.43.

A similar technique will now be applied to a more general class of
functions $f$.

Many special families of orthogonal polynomials have been studied.
In most cases, the knowledge of a special family is considered
satisfactory when an explicit formula for the recurrence coefficients
$a_n$ and $b_n$ as functions of $n$ is associated to a definite
formula for the weight $w$, or measure of orthogonality, see for instance
the final tables of Chihara's book [Chi], whereas the starting point of
the study may be generating functions, Rodrigues formulas, special
functions identities, differential equations, etc.

The simplest way to start the study of the class of {\it semi-classical\/}
orthogonal polynomials is to define them through a differential equation
of their function $f$:

{\bf Definition\/}: The sequence
 $\{ p_n(z)=\gamma_nz^n+\cdots \}_{n=0}^\infty$
is a set of formal semi-classical orthogonal polynomials if
(3) holds  with a function $f$ satisfying the first order linear
differential equation
 $$ Wf' =2Vf +U   \eqno(9) $$
where $W,V$ and $U$ are polynomials ($W\not\equiv 0$).

This is equivalent to the existence of a linear recurrence relation of
the form $\sum_{k=0}^d (n\xi_k+\eta_k)\mu_{n+k}=0$ for the formal
moments $\mu_n$ [BeR].

Moreover, only regular semi-classical orthogonal polynomials will be
considered here, so that $\gamma_n\ne 0, n=0,1,\ldots$

Of course, (4) must be possible with  an expansion of the form (3), so that
degree $U \le \max($ degree $W\ -2,$ degree $V\ -1)$.
All the classical families are recovered
 when degrees of $W$ and $V \le$ 2 and 1.

We will consider especially

{\bf Definition\/}: Generic semi-classical orthogonal polynomials
are semi-classical orthogonal polynomials where $m=$ degree $W  \ge 2$,
degree $V<m$, the zeros $x_1$, $x_2,\ldots x_m$ of $W$ are
distinct, and the residues $\alpha_k=2V(x_k)/W'(x_k)$ are not integers,
$k=1,2,\ldots,m$. The Jacobi polynomials correspond to $m=2$.

We have then:

{\bf Proposition\/}: {\sl Generic semi-classical orthogonal polynomials are
orthogonal with respect to a (possibly complex) {\it generalized
Jacobi weight function\/} \hfill\break
 $w(x)=A_j \Pi_1^m(x-x_k)^{\alpha_k}$ on arcs
$S_j, j=1,2,\ldots m$ of the complex plane.}

Indeed, (9) has exactly one holomorphic solution
 $f_j(z)=c_{j,0}+c_{j,1}(z-x_j)+\cdots$ in a
neighbourhood of the singular point $x_j$, as the equations for the
$c_{j,i}$'s are $2V(x_j)c_{j,0}+U(x_j)=0$ and
$W'(x_j)ic_{j,i}+\cdots = 2V(x_j)c_{j,i}+\cdots, i=1,2,\ldots$ have
exactly one solution, as $V(x_j)\ne 0$ and
$W'(x_j)i-2V(x_j)=W'(x_j)(i-\alpha_j)$ cannot vanish (this can also be
seen as a most elementary application of L.Fuchs theory of linear
differential equations). As $\Pi_1^m(x-x_k)^{\alpha_k}$ is a solution
of the homogeneous equation (9), one has
 $f(z)=f_j(z)+B_j \Pi_1^m(z-x_k)^{\alpha_k}$ near $x_j$, on one side
of the cut. A Cauchy integral expression of $f(z)$ will, after a distorsion
of the integration contour (as in [N] \S1.2), involve the difference
of the  limit functions $f_+$ and $f_-$ which is a multiple
of $\Pi_1^m(x-x_k)^{\alpha_k}$ on a cut. This gives $w$ on $S$.
Let $w(z)$ be a continuation of $w$ on some side of the cut, then we
have
 $$ f(z) = f_j(z) + C_j w(z)  \eqno(10) $$
near $x_j$.

Non generic semi-classical orthogonal polynomials can be considered
as limit cases, for instance, a weight $\exp P(x)$, where $P$ is a
polynomial, is the limit of $(1+P(x)/N)^N$ when $N\to\infty$\ \dots
See [Al] and [Bel] for other proofs and examples.

Anyhow, as $f_+$ and $f_-$ along the two sides of a system of cuts
are solutions of the same equation (9), their difference must be a solution
of the homogeneous equation: {\it semi-classical orthogonal polynomials
are orthogonal with respect to a $($possibly complex$)$ weight function
$w$ satisfying
$$ Ww'=2Vw  \eqno(11) $$
 on a system of cuts\/}, masspoints may also
be present if $f$ has poles. Examples have been given in [BoN], [HvR1],
[HvR2] and [Sho]; the whole class of true positive semi-classical
measures on real sets is given in [BLN].

Conversely, Shohat [Sho] develops the theory starting from a weight
function satisfying (11)  on an interval. Let us generalize this
to a given set of arcs $S$, and show that (9) is recovered: if needed,
we multiply $W$ and $V$ by common factors in order to have
$\lim W(x)w(x)=0$ when $x$ tends to any endpoint (eq. (6) of [Sho]).
Then, from (8), $W(z)f(z)= \int_S W(x)w(x) (z-x)^{-1}dx$ plus a polynomial
($\int_S [(W(z)-W(x))/(z-x)]w(x)\,dx$ is a polynomial in $z$). The
derivative gives
$$ (W(z)f(z))' = -\int_S W(x)w(x) (z-x)^{-2}dx +{\rm\ pol.\ } =
                  \int_S (W(x)w(x))' (z-x)^{-1}dx +{\rm\ pol.\ }, $$
by integration by parts, using $Ww\to 0$ at the endpoints of $S$. As
$(Ww)'=(W'+2V)w$, and $\int_S (W'(x)+2V(x))(z-x)^{-1}dx =
                      (W'(z)+2V(z))\int_S(z-x)^{-1}dx +$ a polynomial,
we find indeed $Wf'=2Vf+$ a polynomial, i.e, (9).

{\bf 3. Differential relations and equations for
formal semi-classical orthogonal polynomials.}

Now, we go further, following Laguerre ([Lag] sec.\ 2, see also [HvR1],
[Per] \S~76): from (4) and (9),
$$ \eqalign{
   0&=W\left[\dfrac{p_{n-1}^{(1)}}{p_n} +\dfrac{\varepsilon_n}{p_n}
      \right]'
   -2V\left[\dfrac{p_{n-1}^{(1)}}{p_n} +\dfrac{\varepsilon_n}{p_n}
      \right] -U  \cr
   &=\dfrac{
      W [p_{n-1}^{(1)}{}'\pph n-p_n'p_{n-1}^{(1)}]
    -2V p_{n-1}^{(1)}\pph n -Up_n^2  }{p_n^2}
   +  W\left[\dfrac{\varepsilon_n}{p_n}  \right]'
   -2V\dfrac{\varepsilon_n}{p_n} \cr
} $$
so,
$$ \Theta_n=    W [p_{n-1}^{(1)}{}'\pph n-p_n'p_{n-1}^{(1)}]
    -2V p_{n-1}^{(1)}\pph n -Up_n^2  \eqno(12) $$
is a polynomial of degree bounded by a constant, as
$$ \Theta_n =
   -p_n^2  W\left[\dfrac{\varepsilon_n}{p_n}  \right]'
   +2V\varepsilon_n  p_n =
  W(\varepsilon_n p_n' - \varepsilon_n' p_n) +2V\varepsilon_n p_n \eqno(13)$$
is bounded by a power $\le\max($ degree $W\ -2$, degree $V\ -1$) for
large argument. For given $W$ and $V$, (13) with (6) and (7) allow
to give $\Theta_n$ in terms of $n$ and the recurrence coefficients
$a$'s and $b$'s. Moreover, expanding (13) up to {\it negative\/}
powers of $z$ yields equations for these coefficients.
This is a first hint towards identities ({\bf Laguerre-Freud's
equations\/})   for the
recurrence coefficients of semi-classical orthogonal polynomials.
See [BeR] for a technique involving Tur\'an determinants.

Identities like (12) involving orthogonal polynomials of arbitrary
high degree on one side and polynomials of bounded degree with respect
to $n$ on the other side occur whenever one has a functional
equation $P(f)=0$ for $f$, provided the elimination of $f$ in
$P(p_{n-1}^{(1)}/p_n +O(z^{-2n-1}))=0$ is simple enough. This happens
if $P$ applied to a rational function $\varphi/\psi$ produces
another rational function with denominator $\xi$ of degree not much larger
that {\it twice\/} the degree of $\psi$. Then, multiplication by
this denominator $\xi$ will produce polynomials and, roughly speaking,
products of $\xi$ and the error term $O(z^{-2n-1})$ which will keep
a small rate of growth at $\infty$. Exemples of valid functionals
$P$ are quadratic polynomials (discussed in the preceding section:
$f=(A-B^{1/2})/C \Rightarrow (Cf-A)^2-B=0$),
linear differential operators of first order discussed here, both
giving $\xi=p_n^2$, and Riccati differential operators (theory of
Laguerre-Hahn orthogonal polynomials [Mag1]). Difference operators may
also be considered, they can leave things like $\xi(z)=p_n(z)p_n(z+h)$,
$\xi(z)=p_n(z)p_n(qz)$, etc. [Mag4]

 In the generic case, let $W(z)=\Pi_1^m (z-x_k) = z^m -(\sum_1^m x_k)
z^{m-1}+\cdots$, then $2V(z)=W(z)\sum_1^m(\alpha_k/(z-x_k)) =
(\sum_1^k \alpha_k)z^{m-1} +[\sum_1^m(\alpha_k x_k) -
(\sum_1^m x_k)(\sum_1^m \alpha_k)]z^{m-2}+\cdots$, using (13), (6)
and (7):
$$\displaylines{
  \Theta_n(z) =
  \left(2n+1+\sum_1^m \alpha_k\right)z^{m-2} + \hfill \cr
+\left[\left(2n+1+\sum_1^m \alpha_k\right)
         \left(b_n-\sum_1^m x_k\right) +2\sum_0^{n-1}b_i +b_n
 +\sum_1^m (\alpha_k x_k)\right] z^{m-3} +\cdots \hfill(14)  \cr
 }
              $$

From (5), replace $\Theta_n$ in (12) by
$(\pph{n-1} p_{n-1}^{(1)}-\pph n p_{n-2}^{(1)})a_n\Theta_n$:
\hfill\break
$p_{n-1}^{(1)}[W\pph n{}'+V\pph n +a_n\Theta_n\pph{n-1}]=
 \pph n{}[Wp_{n-1}^{(1)}{}'-Vp_{n-1}^{(1)}+a_n\Theta_n p_{n-2}^{(1)}
-U\pph n ],$
which must therefore have the form $\Omega_n \pph n p_{n-1}^{(1)}$,
where $\Omega_n$ is a new auxiliary polynomial of bounded degree.
Using again (5), one has
$$\Omega_n = a_nW[p_{n-1}^{(1)}{}' \pph{n-1}-p_n' p_{n-2}^{(1)}]
 -a_nV[p_{n-1}^{(1)}\pph{n-1}+\pph n p_{n-2}^{(1)}]
 -a_nU\pph n \pph{n-1} $$
And, with (4):
$$\Omega_n = a_nW(\varepsilon_{n-1}p_n'-\varepsilon_n'p_{n-1})
            +a_nV(\varepsilon_{n-1}p_n +\varepsilon_n p_{n-1}).
  \eqno(15) $$
This yields the two {\it differential relations\/}:
$$\eqalign{
  Wp_n' &= (\Omega_n-V)p_n -a_n\Theta_n p_{n-1} \cr
 Wp_{n-1}^{(1)}{}' &= (\Omega_n+V)p_{n-1}^{(1)}
                   -a_n\Theta_n p_{n-2}^{(1)} +U\pph n\cr
}$$
We get rid of the $Up_n$ term of the second equation by forming an
equation for $fp_n$, using (9), and subtracting the second equation:
$ W\varepsilon_n'  = (\Omega_n+V)\varepsilon_n -a_n\Theta_n
 \varepsilon_{n-1}$. We recover the form of the first equation by
using (11):
$ W(\varepsilon_n/w)'  = (\Omega_n-V)\varepsilon_n/w -a_n\Theta_n
 \varepsilon_{n-1}/w$. In order to have a differential system, we have
to give $y_{n-1}'$ ($y=p$ or $\varepsilon/w$) in terms of $y_n$ and
$y_{n-1}$. As $y_n$  satisfies the recurrence relations (2),
$ Wy_{n-1}'  = (\Omega_{n-1}-V)y_{n-1} -a_{n-1}\Theta_{n-1}y_{n-2}$
turns easily as
$ Wy_{n-1}'  = a_n\Theta_{n-1} y_n +
 (\Omega_{n-1}-V-(z-b_{n-1})\Theta_{n-1})y_{n-1} $.
From (15), (13) and (2),
 $$\Omega_{n+1}(z) = (z-b_n)\Theta_n(z)-\Omega_n(z), \eqno(16)$$
so we finally have the {\it differential system\/}:

$$ Y'=AY : \quad
  \left[\matrix{ p_n & \varepsilon_n/w \cr
                 p_{n-1} & \varepsilon_{n-1}/w \cr
               }\right]' =
     \dfrac 1W\,
  \left[\matrix{ \Omega_n-V & -a_n\,\Theta_n  \cr
                 a_n\Theta_{n-1} & -\Omega_n\,-V  \cr
               }\right]
  \left[\matrix{ p_n & \varepsilon_n/w \cr
                 p_{n-1} & \varepsilon_{n-1}/w \cr
               }\right]\quad .
 \eqno(17) $$

This differential system gives the whole differential history of the
semi-classical orthogonal polynomials. Laguerre [Lag] and many other
people ([AtE] [Ha1] [Ha2] [Nev] [Sho] etc.~)
have preferred the scalar second order form obtained from
eliminating $y_{n-1}$ in $Wy_n'=(\Omega_n-V)y_n-a_n\Theta_n y_{n-1}$
and $Wy_{n-1}'=a_n\Theta_{n-1} y_n -(\Omega_n+V)y_{n-1}$:
$$W\Theta_n y_n'' = (W\Theta_n' -W'\Theta_n -2V\Theta_n)y_n' +K_ny_n,
  \eqno(18) $$
with $K_n=(\Omega_n-V)'\Theta_n-(\Omega_n-V)\Theta_n'
 +\Theta_n(\Omega_n^2-V^2-a_n^2\Theta_n \Theta_{n-1})/W$,
which is a polynomial, as putting
$a_{n+1}y_{n+1}'=(z-b_n)y_n'+y_n-a_ny_{n-1}'$  (derivative of (2)) in
$a_{n+1}Wy_{n+1}'=a_{n+1}(\Omega_{n+1}-V)y_{n+1}-a_{n+1}^2\Theta_{n+1}
y_n$, using again (2), and the differential equation (17) for
$Wy_{n-1}'$ gives an expression of the form
$Ay_n=By_{n-1}$, with $B=0$ from (16), whence $A=0$, which is
$$ (z-b_n)(\Omega_{n+1}-\Omega_n) = W +a_{n+1}^2\Theta_{n+1}-a_n^2
 \Theta_{n-1}. \eqno(19)$$
Multiplying by (16) and summing on $n$, one finds
 $$\Omega_n^2 -a_n^2\Theta_n\Theta_{n-1} = V^2+W\sum_0^{n-1}\Theta_i ,
\eqno(20)$$
knowing that $\Omega_0=V$.

With $z_n=(Ww/\Theta_n)^{1/2}y_n$, we have a form without first derivative
$$z_n''=\left\{ \dfrac34\left(\dfrac{\Theta_n'}{\Theta_n}\right)^2
               -\dfrac12      \dfrac{\Theta_n''}{\Theta_n}
               -\dfrac12      \dfrac{\Theta_n'}{\Theta_n}
                             \dfrac{W'+2\Omega_n}W
                         +   \dfrac{4V^2-{W'}^2}{4W^2}
                         +   \dfrac{W''+2\Omega_n'}{2W}
                        +\dfrac{\sum_0^{n-1}\Theta_i}W
         \right\}z_n,  \eqno(21)$$
used by R.\ Fuchs [RFu] in the case $m=$ degree $W=3$.

Laguerre ([Lag], see also [GaN]) finds equations for the recurrence
coefficients and the coefficients of $\Theta_n$ and $\Omega_n$ by
using (16) and (19), keeping the degrees of $\Theta_n$ and
$\Omega_n$ bounded when $n$ increases. We may express everything in
terms of the recurrence coefficients alone, then the expansion
of $\Omega_n$, constructed on the same lines as (14), will be
useful:
$$\displaylines{
  \Omega_n(z) =
   \left[n+(\sum\nolimits_1^m \alpha_k)/2\right] z^{m-1} + \hfill\cr
+\left[\sum_0^{n-1}b_i -n\sum\nolimits_1^m x_k +
                    \left( \sum\nolimits_1^m(\alpha_k x_k) -
    (\sum\nolimits_1^m x_k)(\sum\nolimits_1^m \alpha_k)\right)/2\right]
     z^{m-2}+     \hfill          \cr
+\left[\sum_0^{n-1}b_i^2 +2\sum_1^{n-1}a_i^2 -(\sum\nolimits_1^m x_k)
  \left(\sum_0^{n-1}b_i + \sum\nolimits_1^m(\alpha_k x_k)/2\right)
        +\right.\hfill\cr
        \left.
 +(n+\sum\nolimits_1^m \alpha_k/2)\left(\sum_{k<\ell\le m}x_k x_\ell
    \right) +\sum\nolimits_1^m \alpha_k x_k^2/2 +
   \left(2n+1+\sum\nolimits_1^m \alpha_k\right)a_n^2 \right]z^{m-3} +\ldots
     \hfill(22) \cr
           }
              $$

Consider for instance the case $m=3$ (simplest generalized Jacobi
polynomials): from (14), $\Theta_n$ is a polynomial of degree 1 with
a known coefficient of $z$ and a constant coefficient depending on
the $b_i$'s up to $b_n$; from (22), $\Omega_n$ is a polynomial of
degree 2 with a known coefficient of $z^2$ and two other coefficients
depending on the $b_i$'s and the $a_i$'s up to the index $n-1$ (see
example 1 in section 5) . The
constant coefficients of (16) and (20) give nonlinear relations
for $a_n$ and $b_n$. The meaning of the solutions of these recurrence
relations for the recurrence coefficients of (2) is not obvious. Even
the simplest relations found in nongeneric cases (as (30)
or (31)) are baffling.

The explanation in terms of Painlev\'e transcendents and similar
functions, i.e., solutions of remarkable high-order nonlinear
differential equations in terms of a well-chosen parameter, will
be given now. The derivation is based on the isomonodromy properties
of (18). Later on, examples will show that a more elementary
derivation is possible.

{\bf 4. Monodromy matrices and isomonodromy identities.}

D.\ \&
G.\ Chudnovsky remarked ([Chu2], see also (5.1.18) in [N])  how (18) has a form
already
investigated in the period 1890-1910 by authors working on
isomonodromy deformations ([RFu], [Pain]).

Let $Y(z)$ be a fundamental matrix of solutions of the differential
system $Y'(z)=A(z)Y(z)$, defined outside a system of cuts joining
the singular points (poles of $A$) of the equation. When $z$ follows
a contour about a singular point $x_j$, let us solve $Z'(z)=A(z)Z(z)$
with the initial value $Z(z_0)=Y(z_0)$ at a starting point on the contour.
As long as no cut is crossed, $Z(z)=Y(z)$. This is no more true when
one or several cuts are crossed but, when we come back in a neighbourhood
of $z_0$, the columns of the matrix of solutions $Z(z)$ must be fixed
combinations of the columns of the initial fundamental matrix of
solutions: $Z(z)=Y(z)M_j$. This matrix $M_j$ is called the
{\bf monodromy matrix\/} of $Y'=AY$ at the singular point $x_j$ (only
{\it regular\/} singularities are considered here).

{\bf Theorem~1.\/} {\sl Generic formal semi-classical orthogonal polynomials
satisfy differential systems (17) with monodromy matrices
 $$ M_j = \left[ \matrix {1  &  C_j[1-\exp(-2\pi i\alpha_j)] \cr
                          0  & \exp(-2\pi i\alpha_j)  \cr
                         }\right]
 $$
at the singular points $x_j, j=1,2,\ldots m$.}

Indeed, $p_n$ and $p_{n-1}$ are not modified after a circle about
$x_j$, but $\varepsilon_n$, $\varepsilon_{n-1}$ and $w$ have a
branchpoint there. According to the discussion made in the proof of
(10), $f(z)=f_j(z)+B_j\Pi_1^m(z-x_k)^{\alpha_k}$ with some
determination of the powers near $x_j$, near a side of a cut.
By following a contour about $x_j$, $f_j$ returns to its previous value,
but    $\Pi_1^m(z-x_k)^{\alpha_k}$ has been multiplied by
$\exp(2\pi i \alpha_j)$. The same happens with $w$. Therefore,
from (10),
$\varepsilon_n/w=(f\pph n-p_{n-1}^{(1)})/w =
                 (f_j\pph n-p_{n-1}^{(1)})/w+C_j \pph n$  becomes
    $\exp(-2\pi i\alpha_j)(f_j\pph n-p_{n-1}^{(1)})/w+C_j \pph n$ =
  $\exp(-2\pi i\alpha_j)\varepsilon_n/w + [1-\exp(-2\pi i \alpha_j)]
  C_j p_n$.

{\it
This shows that the monodromy matrices of (17) at the singular
points remain unchanged if the exponents $\alpha_k$ remain
unchanged and if the weight $w$ on $S$ is adapted so that the multipliers
$C_k$ remain the same.\/} However, one may vary the positions of the singular
points $x_k$. The quantities $f$, $p_n$, $a_n$, $b_n$, $\Theta_n$ etc.
will then be subject to extrememely interesting {\it isomonodromy
deformations\/}. Here is a sketch ([LD] III, from p.128 onwards),
applied to the specific equation (17):

Let the $x_k$ depend on a single parameter $t$, and let us define the
matrix
$$ H= \dfrac{\partial Y}{\partial t}\, Y^{-1} ,       $$
as $\partial M_j/\partial t=0$, $H$ does not change when $z$ achieves
a contour about $x_j$. So, $H$ has no branchpoints at the $x_j$'s.
To get a better view of what happens at the singular points, we expand
$H$ (using $\det Y = 1/(a_n w)$, from (5)):

$$ H=a_n\left[\matrix{
     \dot p_n \varepsilon_{n-1}-p_{n-1}\dot \varepsilon_n
       +p_{n-1}\varepsilon_n \dot{w}/w  &
    -\dot p_n \varepsilon_{n}+p_{n}\dot \varepsilon_n
       -p_{n}\varepsilon_n \dot{w}/w  \cr
     \dot p_{n-1} \varepsilon_{n-1}-p_{n-1}\dot \varepsilon_{n-1}
       +p_{n-1}\varepsilon_{n-1} \dot{w}/w  &
    -\dot p_{n-1} \varepsilon_{n}+p_{n}\dot \varepsilon_{n-1}
       -p_{n}\varepsilon_{n-1} \dot{w}/w  \cr
                   } \right],        \eqno(23)
$$
where the dot derivative is $\partial /\partial t$.
From (4) and (10), one has $\varepsilon_n = \varepsilon_{n,j}
  + C_j w p_n$ near $x_j$, where $\varepsilon_{n,j}$ is regular
near $x_j$. The singular terms cancel nicely in the combinations
of (23)  (remember that $\dot C_j=0!$);  the ratio $\dot w/w$ has a simple
pole at $x_j$ with residue $-\alpha_j \dot x_j$ (as $\dot\alpha_j=0$).
We are left with
$$ H = H_\infty + \sum_{j=1}^m  H_j (z-x_j)^{-1}, $$
with
 $$ H_j = -\alpha_j \dot x_j a_n \left[ \matrix{
            p_{n-1}\varepsilon_{n,j} & -p_n\varepsilon_{n,j} \cr
            p_{n-1}\varepsilon_{n-1,j} & -p_n\varepsilon_{n-1,j} \cr
                 } \right] \qquad j=1,\ldots, m,
$$
where the $ p_r \varepsilon_{s,j}$'s are the values at $z=x_j$. As
$W(x_j)=0$, (13) tells that $\Theta_n=2V\varepsilon_{n,j}p_n$ at
$z=x_j$, and (15) with (5) gives
$\Omega_n = V+2a_nV\varepsilon_{n,j}p_{n-1} =
           -V+2a_nV\varepsilon_{n-1,j}p_n$ at $x_j$. With
$\alpha_j = 2V/W'$ at $x_j$, one finds from (17):
$$ A=\sum_{j=1}^m (z-x_j)^{-1} A_j \qquad\Rightarrow\qquad
   H=H_\infty -\sum_{j=1}^m (z-x_j)^{-1} \dot x_j A_j,$$
A direct inspection of (23) when $z\to\infty$ gives,
using (6) and (7),
$$ H_\infty = \left[\matrix{ \dot\gamma_n/\gamma_n & 0 \cr
                      0 &   -\dot\gamma_{n-1}/\gamma_{n-1}  \cr
                           }\right]
$$
in the generic case, as $\dot w/w  =   -\sum_1^m \alpha_k \dot x_k
/(z-x_k)\to 0$ when $z\to\infty$.

Finally, the {\it differential equations in $t$\/} appear by working
$$ \eqalign{
   \partial^2 Y/\partial z\partial t &=\partial \dot Y/\partial z =
         (HY)' = H'Y+HY' = H'Y+HAY = \cr
 = \partial^2 Y/\partial t\partial z &=\partial Y' /\partial t =
     \dot{(AY)} = \dot A Y + A\dot Y = \dot A Y +AHY , \cr
}$$
whence
$$ \dot A = H'+HA-AH. \eqno(24) $$

\rightline{\sl C'\'etaient les cieux ouverts}
\rightline{    Stendhal\ \ \ \ \ \ \ \ \ \ \ }

This equation (24) has an incredibly inspiring form, explaining how
this theory is related to integrable Hamiltonians, B\"acklund
transformations, Lax pairs, Toda lattices, solitons, etc.
whereas the connection with orthogonal
polynomials, special functions,  continued fractions, Diophantine
approximations has been worked with great virtuosity by G.\ \& D.\
Chudnovsky [Chua] [Chub] [Chu0] [Chu1] [Chu2],\dots

In the generic case, we have for the residue matrices
$$ \dot A_j = H_\infty A_j- A_j H_\infty +
  \sum_{ \scriptstyle{k=1}\atop\scriptstyle{k\ne j} }^{k=m}
    \dfrac{\dot x_j -\dot x_k}{x_j-x_k} \, (A_k A_j - A_j A_k),
 \qquad j=1,\ldots,m $$
called the {\it Schlesinger equations\/} (see [Chua]).

Now, we show how these equations lead to differential equations
for the recurrence coefficients:

{\bf Theorem~2.\/} {\sl Let $a_n$ and $b_n$ be the recurrence coefficients
of $(2)$ for generalized Jacobi orthogonal polynomials related to a
$($ possibly complex $)$ weight of the form
$\Pi_1^m (x-x_j)^{\alpha_j}$, on a set of arcs joining the $x_j$'s,
 where at least one of the $x_j$'s depend
on a parameter $t$. Then, we have the {\it Toda equations\/}
$$ \dfrac{\dot a_n}{a_n} = \dfrac12 \sum_{k=1}^m
    \dfrac{(\Theta_n(x_k)-\Theta_{n-1}(x_k))\dot x_k}{W'(x_k)},
  \eqno(25) $$
$$ \dot b_n       =          \sum_{k=1}^m
    \dfrac{(\Omega_{n+1}(x_k)-\Omega_n(x_k))\dot x_k}{W'(x_k)},
  \eqno(26) $$
where $W(x)=\Pi_1^m (x-x_k)$, and $\Theta_n$ and $\Omega_n$
are polynomials introduced in $(12)- (15)$. }

Indeed, from (17), the residue matrix $A_j$ is
$$ A_j = \dfrac1{W'(x_j)} \left[ \matrix{
         \Omega_n(x_j)-V(x_j)  & -a_n\Theta_n(x_j) \cr
        a_n\Theta_{n-1}(x_j)  &  -\Omega_n(x_j)-V(x_j) \cr}\right] $$
we have
$$ H_\infty A_j -A_j H_\infty =
 - \dfrac{\dot\gamma_n/\gamma_n + \dot\gamma_{n-1}/\gamma_{n-1}}
    {W'(x_j)}\, a_n \left[\matrix{ 0 & \Theta_n(x_j) \cr
                               \Theta_{n-1}(x_j) & 0 \cr } \right],$$
$$\displaylines{
  A_k A_j - A_j A_k = \dfrac{a_n}{W'(x_j)W'(x_k)} \ \times \hfill\cr
 \times  \left[\matrix{
 a_n(\Theta_n(x_j)\Theta_{n-1}(x_k)-\Theta_n(x_k)\Theta_{n-1}(x_j)) &
2(\Theta_n(x_k)\Omega_n(x_j)-\Theta_n(x_j)\Omega_n(x_k)) \cr
2(\Theta_{n-1}(x_k)\Omega_n(x_j)-\Theta_{n-1}(x_j)\Omega_n(x_k)) &
 a_n(\Theta_n(x_k)\Theta_{n-1}(x_j)-\Theta_n(x_j)\Theta_{n-1}(x_k)) \cr
  }\right].\cr}$$
The Schlesinger equations for the off-diagonal elements of $A_j$ are
$$ -\dfrac\partial{\partial t} \dfrac{\Theta_n(x_j)}{W'(x_j)} =
 -2\dfrac{\dot\gamma_n}{\gamma_n}  \dfrac{\Theta_n(x_j)}{W'(x_j)} +2
\sum_{k\ne j} \dfrac{\dot x_j -\dot x_k}{x_j-x_k}
    \dfrac{\Theta_n(x_k)\Omega_n(x_j)-\Theta_n(x_j)\Omega_n(x_k)}
      {W'(x_j)W'(x_k)},\eqno(27)$$
$$ \dfrac\partial{\partial t} \dfrac{\Theta_{n-1}(x_j)}{W'(x_j)} =
 -2\dfrac{\dot\gamma_{n-1}}{\gamma_{n-1}}  \dfrac{\Theta_{n-1}(x_j)}{W'(x_j)}
  +2
\sum_{k\ne j} \dfrac{\dot x_j -\dot x_k}{x_j-x_k}
 \dfrac{\Theta_{n-1}(x_k)\Omega_n(x_j)-\Theta_{n-1}(x_j)\Omega_n(x_k)}
      {W'(x_j)W'(x_k)},$$
where $a_n\gamma_n=\gamma_{n-1} \Rightarrow \dot a_n/a_n
 +\dot\gamma_n/\gamma_n = \dot\gamma_{n-1}/\gamma_{n-1}$ has been
used. Increasing $n$ by 1 in the second equation and adding to the
first one,
$$ 0= -4\dfrac{\dot\gamma_n}{\gamma_n} \dfrac{\Theta_n(x_j)}{W'(x_j)}
 +2\sum_{k\ne j} (\dot x_j-\dot x_k)
   \dfrac{\Theta_n(x_j)\Theta_n(x_k)}{W'(x_j)W'(x_k)},$$
where (16) has been used. At this point, we don't have to avoid
the term $k=j$ anymore in the sum. Moreover, as any polynomial
$P(z)=\pi_0 z^{m-1}+\cdots$ satisfies $\pi_0=\sum_1^m P(x_k)/W'(x_k)$,
(coefficient of $z^{-1}$ in $P(z)/W(z) = \sum_1^m P(x_k)/((W'(x_k)(z-x_k))$),
and as the degree of $\Theta_n$ is $m-2$ ((14)), $\dot x_j$
disappears from the sum:
$$ \dfrac{\dot\gamma_n}{\gamma_n} = -\dfrac12 \sum_{k=1}^m
    \dfrac{\Theta_n(x_k)\,\dot x_k}{W'(x_k)}  $$
and (25) follows from $a_n\gamma_n=\gamma_{n-1}$.\hfill\break
Now, we come to the first diagonal element of the Schlesinger's
equations:
$$\dfrac\partial{\partial t} \dfrac{\Omega_n(x_j)-V(x_j)}
    {W'(x_j)} =
a_n^2 \sum_{k\ne j} \dfrac{\dot x_j -\dot x_k}{x_j-x_k}
 \dfrac{\Theta_n(x_j)\Theta_{n-1}(x_k)-\Theta_n(x_k)\Theta_{n-1}(x_j)}
      {W'(x_j)W'(x_k)},\eqno(28)$$
for $j=1,\ldots, m$. As
$\dfrac{\Theta_n(x)\Theta_{n-1}(y)-\Theta_n(y)\Theta_{n-1}(x)}
 {x-y}$ is some polynomial, say $\sum_{p,q} \tau_{p,q}x^p y^q$, of
degree $\max(p,q) < m-2$ in $x$ and $y$, the term $k=j$ may be
included in the sum as before. Still using
 $\sum_1^m P(x_k)/W'(x_k)=0$ for polynomials $P$ of degree less than
$m-1$, $\dot x_j$ may also be removed, and
$$\eqalign{
 \sum_{j=1}^m \dfrac1{z-x_j} \dfrac\partial{\partial t}
   \dfrac{\Omega_n(x_j)-V(x_j)}{W'(x_j)} &=
   -a_n^2 \sum_{k=1}^m \dfrac{\dot x_k}{W'(x_k)}
  \sum_{p,q} \tau_{p,q} x_k^q \sum_{j=1}^m
    \dfrac{x_j^p}{(z-x_j)W'(x_j)} \cr
  &=   -a_n^2 \sum_{k=1}^m \dfrac{\dot x_k}{W'(x_k)}
  \sum_{p,q} \tau_{p,q} x_k^q \dfrac{z^p}{W(z)} \cr
  &=   -a_n^2 \sum_{k=1}^m \dfrac{\dot x_k}{W'(x_k)}
 \dfrac{\Theta_n(z)\Theta_{n-1}(x_k)-\Theta_n(x_k)\Theta_{n-1}(z)}
       {W(z)\,(z-x_k)}.\cr
          }
$$
The left-hand side is
$$\displaylines{
   \dfrac\partial{\partial t}\left(
  \sum_{j=1}^m\dfrac1{z-x_j}  \dfrac{\Omega_n(x_j)-V(x_j)}{W'(x_j)} \right)
 -
 \sum_{j=1}^m \dfrac\partial{\partial t}\left(\dfrac1{z-x_j}\right)
   \dfrac{\Omega_n(x_j)-V(x_j)}{W'(x_j} = \hfill\cr
 \hfill =
   \dfrac\partial{\partial t}\, \dfrac{\Omega_n(z)-V(z)}{W(z)}  -
 \sum_{k=1}^m \dfrac{\dot x_k}{(z-x_k)^2}\,
   \dfrac{\Omega_n(x_k)-V(x_k)}{W'(x_k}, \cr
} $$
and we take the $z^{-2}$ term in the expansion about $\infty$,
 the right-hand
side vanishes as degree $\Theta_n < m-1$, and using (22) in the
left-hand side:
$$ \sum_0^{n-1}\dot b_i  - \sum_1^m \dot x_k \dfrac{\Omega_n(x_k)-V(x_k)}
   {W'(x_k)} =0   $$
yields (26).

In concrete situations, (25) and (26) will be used, together with
other non differential identities (Freud Laguerre equations for the
recurrence coefficients), but we may prefer to return to (24), or
even use {\it ad hoc\/} differential relations. In the generic case,
(27) and (28) for $j=1,2,\ldots,m$ give a differential system for
$2m$ unknowns $\Theta_n(x_j)$ and $\Omega_n(x_j)$, $j=1,2,\ldots,m$,
when the quantities $a_n^2\Theta_{n-1}(x_j)$ are eliminated with
the help of (20) at $x_j$ (recall that $W(x_j)=0$). However, considering
from (14) and (22) that there are only $2m-3$ unknown coefficients
in $\Theta_n$ and $\Omega_n$, further eliminations are possible.
We start with an example of generic semi-classical weight with $m=3$.

\bigskip

{\bf 5. Example~1. Generalized Jacobi weight with three factors\hfill\break
  $(1-x)^\alpha x^\beta (t-x)^\gamma$.}

So, $W(z)=z(z-1)(z-t)$, $V(z)=(\alpha z(z-t)+\beta(z-1)(z-t)+
\gamma z(z-1))/2$, the support $S$ joins $0,1,$ and $t$ in some way,
or is an arc joining only two of these points.
(14) and (22) yield readily
$$ \Theta_n(z)=\nu_n z+\vartheta_n,\qquad \Omega_n(z)=
      \dfrac{\nu_n -1}2 z^2
   +\kappa_n z+\omega_n,$$
with $\nu_n=2n+1+\alpha+\beta+\gamma$,
$\vartheta_n = \nu_n(b_n-1-t)+2\sum_0^{n-1}b_i +b_n+\alpha+\gamma t$,
$\kappa_n = \sum_0^{n-1}b_i -(\nu_n -1)(1+t)/2 +(\alpha+\gamma t)/2$,
and
$\omega_n = \sum_0^{n-1}(b_i^2 -(t+1)b_i +2a_i^2)-(t+1)(\alpha+\gamma t)/2
+(\nu_n-1)t/2 +(\alpha+\gamma t^2)/2 +\nu_n a_n^2$.

\noindent
(25) and (26) are here, with $\dot x_k=\delta_{k,3}$,
   $W'(x_3)=W'(t)=t(t-1)$,
 $$ \dfrac{\dot a_n}{a_n} = \dfrac{-2+(\nu_n+1)b_n-(\nu_n-3)b_{n-1}}
                                  {2t(t-1)},\
    \dot b_n  = \dfrac{b_n(b_n-1)+(\nu_n+2)a_{n+1}^2-(\nu_n-2)a_n^2}
                                  {t(t-1)}.\
$$
One would have a true differential system if $b_{n-1}$ and
$a_{n+1}^2$ were simple functions of $b_n$ and $a_n$, but this does not
seem to be the case here. So, we try with the unknowns $\vartheta_n$,
$\kappa_n$ and $\omega_n$ instead. In (27), using
$[\Theta_n(x)\Omega_n(y)-\Theta_n(y)\Omega_n(x)]/(y-x)=
    (\nu_n -1)\nu_n xy/2+\vartheta_n(\nu_n-1)(x+y)/2+\zeta_n$,
with $\zeta_n=\vartheta_n \kappa_n - \nu_n \omega_n$, with $x_j=0,1,t$,
one finds three equations which are all equivalent to
$$ \dot\vartheta_n = \dfrac{-\vartheta_n-\vartheta_n^2 +2\zeta_n}
      {t(t-1)}.\eqno(29)$$

\noindent
In (28), using
$[\Theta_n(x)\Theta_{n-1}(y)-\Theta_n(y)\Theta_{n-1}(x)]/(x-y)=
  \nu_n\vartheta_{n-1}-\vartheta_n\nu_{n-1}$, one finds
two independent equations
$$\eqalign{
   \dot\omega_n &= \dfrac{\omega_n}t -
        \dfrac{a_n^2(\nu_n\vartheta_{n-1}-\vartheta_n\nu_{n-1})}{t(t-1)},\cr
   \dot\kappa_n &= \dfrac{\nu_n-1}{2(t-1)}+\dfrac{\kappa_n}{t-1} +
       \dfrac{\omega_n}{t(t-1)}.\cr
}$$
Now, the three non differential equations (20) at $x=0,1,t$ allow
to eliminate $a_n^2\vartheta_{n-1}$, $a_n^2\nu_{n-1}$:
$$a_n^2\vartheta_{n-1}=\dfrac{\omega_n^2-\beta^2t^2/4}{\vartheta_n},
 a_n^2\nu_{n-1}=\dfrac{((\nu_n-1)/2+\kappa_n+\omega_n)^2
       -\alpha^2(t-1)^2/4}{\nu_n+\vartheta_n}-a_n^2\vartheta_{n-1},$$
and a third equation allowing to eliminate either $\kappa_n$ or
$\omega_n$, actually it is simpler to give everything in function of
$\zeta_n$:
from $a_{n-1}^2[\vartheta_{n-1}/t-(\vartheta_{n-1}+\nu_{n-1})/(t-1)+
  (\vartheta_{n-1}+\nu_{n-1} t)/(t(t-1))]=0$,
$$\displaylines{
  \omega_n = -\dfrac{\alpha^2\vartheta_n (t-1)/4}{(\nu_n-1)(\nu_n+\vartheta_n)}
             +\dfrac{\beta^2 t/4}{\nu_n-1}
        +\dfrac{\gamma^2\vartheta_n t(t-1)/4}{(\nu_n-1)(\nu_n t+\vartheta_n)}-
\hfill\cr\hfill
 -\dfrac{(\nu_n-1)\vartheta_n[\nu_n t(t+1)+\vartheta_n(t^2+t+1)]/4
             +[\nu_n t+\vartheta_n(t+1)]\zeta_n +\zeta_n^2/(\nu_n-1)}
         {(\nu_n+\vartheta_n)(\nu_n t+\vartheta_n)}.  \cr}$$
This allows to give $\dot\zeta_n$ as a function of $\vartheta_n$ and
$\zeta_n$, so to complete (29):
$$\hskip-15pt\eqalign{
   \dot\zeta_n &= \dot\vartheta_n \kappa_n + \vartheta_n \dot\kappa_n
         -\nu_n\dot\omega_n \cr
 &= \dfrac{-\vartheta_n -\vartheta_n^2 +2\zeta_n}{t(t-1)}\, \kappa_n
   +  \dfrac{(\nu_n-1)\vartheta_n}{2(t-1)}+
   \dfrac{\kappa_n\vartheta_n}{t-1} +
       \dfrac{\omega_n\vartheta_n}{t(t-1)}
   - \dfrac{\nu_n\omega_n}t +
 \dfrac{a_n^2\nu_n(\nu_n\vartheta_{n-1}-\vartheta_n\nu_{n-1})}{t(t-1)}.\cr
    }$$
Using the preceding calculations, $a_n^2\vartheta_{n-1}$ and
$a_n^2\nu_{n-1}$ are replaced in terms of $\kappa_n$ and $\omega_n$,
then $\kappa_n=(\zeta_n+\nu_n\omega_n)/\vartheta_n$ is used, and
$\omega_n$ is finally replaced as a function of $\zeta_n$, and what comes
out is
$$\displaylines{\hskip-1pt
\dot\zeta_n = \dfrac1{t(t-1)}\left\{
  \dfrac{\alpha^2(t-1)\vartheta_n(\nu_n t+\vartheta_n)}{4(\nu_n+\vartheta_n)}
 -\dfrac{\beta^2 t(\nu_n+\vartheta_n)(\nu_n t+\vartheta_n)}{4\vartheta_n}
  + \right.\hfill\cr  \qquad
 +\dfrac{(1-\gamma^2)t(t-1)\vartheta_n(\nu_n +\vartheta_n)}
                                                {4(\nu_n t+\vartheta_n)}
 +\left(\dfrac1{\nu_n +\vartheta_n}+\dfrac1{\vartheta_n}
       +\dfrac1{\nu_n t +\vartheta_n}\right)
                 (\zeta_n-\vartheta_n(\vartheta_n+1)/2)^2 +\hfill\cr
    +\left(2\vartheta_n+1+\dfrac{\nu_n t(t-1)}{\nu_n t+\vartheta_n}\right)
                 (\zeta_n-\vartheta_n(\vartheta_n+1)/2)
   +\left. \dfrac{\vartheta_n(\nu_n+\vartheta_n)(\nu_n t+\vartheta_n)}
                 {4t(t-1)}
    {\vphantom{\dfrac{\alpha^2}4}}\right\} \cr
}
$$

\noindent
whence, at last, with
$\zeta_n-\vartheta_n(\vartheta_n+1)/2=t(t-1)\dot\vartheta_n/2$:
$$\displaylines{
  \ddot\vartheta_n = \dfrac1{t(t-1)}\,(-2t\dot\vartheta_n
           -2\vartheta_n\dot\vartheta_n +2\dot\zeta_n) =\hfill\cr
\ \ \ =
  \dfrac12\left(\dfrac1{\nu_n +\vartheta_n}+\dfrac1{\vartheta_n}
       +\dfrac1{\nu_n t +\vartheta_n}\right)\,\dot\vartheta_n^2
   -\left(\dfrac1t+\dfrac1{t-1}-\dfrac{\nu_n}{\nu_n t+\vartheta_n}\right)
   \dot\vartheta_n +\hfill\cr
  \ \ \ +
  \dfrac{\alpha^2\vartheta_n(\nu_n t+\vartheta_n)}{2t^2(t-1)(\nu_n+\vartheta_n)}
 -\dfrac{\beta^2 (\nu_n+\vartheta_n)(\nu_n t+\vartheta_n)}{2t(t-1)^2\vartheta_n}
 +\dfrac{(1-\gamma^2)\vartheta_n(\nu_n +\vartheta_n)}{2t(t-1)
                                                (\nu_n t+\vartheta_n)}
 +\hfill\cr        \ \ \
   + \dfrac{\vartheta_n(\nu_n+\vartheta_n)(\nu_n t+\vartheta_n)}
                 {2t^2(t-1)^2}. \hfill\cr
}
$$
which is a Painlev\'e equation of the sixth kind ([In] \S~14.4)
in $-\vartheta_n/\nu_n$ (the
zero of $\Theta_n$) ([Chua] p.399-402, explaining works of R.Fuchs
on equations of form (21)).

We can return to $a_n$ and $b_n$ as functions of $\vartheta_n$ and
$\zeta_n$ by using again $a_n^2\nu_{n-1}$ as a function of  $\vartheta_n$,
$\kappa_n$ and $\omega_n$ (and $\nu_{n-1}=2n+\alpha+\beta+\gamma-1$ is
known) and taking $b_n$ from $2\kappa_n-\vartheta_n=-(2\nu_n -1)(1+t)
-(\nu_n +1)b_n$. Inverting the connection should give a (probably
algebraic) differential system involving only $a_n$ and $b_n$
(will somebody do that?)

{\bf 6. Example~2. $\exp(x^3/3+tx)$ on $\{x: x^3<0\}$.}

Much simpler identities occur when the weight $w$ is the
exponential of a polynomial, so that $w'/w$ is a polynomial
itself. Recall (end of Section~2) that $W(x)w(x)\to 0$ when $x$
tends to the endpoints (if any) of the support $S$.
We want the simplest case ($W(x)=1$), so that the support cannot
have finite endpoints, but must end on directions where $w(x)\to 0$,
with at least one complex direction (or else all the moments vanish).
So, we can take the set $\{x: x^3<0\}$, or for instance only
$\{x: \arg x=\pm 2\pi/3\}$, or also some equivalent contour, as
$\{x:$ Re $x$ = a positive constant$\}$ leading to Airy functions and
integrals ([Chu2] [Mar1]).
The weight can be considered as a confluent generalized Jacobi
weight with singular points at $\infty$:
$w(z)=\lim_{N\to\infty} [1+(z^3/3+tz)/N]^N$, with an exponent $N$
independent of the parameter $t$.  As (25) and (26) hold for
distinct finite singular points, we return to (24) assumed to be
still valid:
here, $W(z)=1$, $2V(z)=w'(z)/w(z)=z^2+t$. Working (13) and (15)
about $\infty$, we have
$$ \Theta_n(z)=z+b_n,\qquad \Omega_n(z)=(z^2+t)/2+a_n^2. $$
Pushing (13) and (15) up to the $z^{-1}$ term, one finds the
corresponding Laguerre-Freud equations, i.e., the identities
$$a_n^2+a_{n+1}^2+b_n^2+t=0,\qquad n+a_n^2(b_n+b_{n-1})=0.\eqno(30)$$
We compute $H$ in (23) up to the $O(1)$ term, as $H$ is now expected
to be a polynomial (see  [Fed] \S~2), taking care of
$\dot w/w=z$:
$$A= \left[\matrix{ a_n^2         &     -a_n(z+b_n)     \cr
                  a_n(z+b_{n-1}) &   -a_n^2-z^2-t \cr  }\right],\
  H= \left[\matrix{\dot\gamma_n/\gamma_n  &     -a_n      \cr
                  a_n &  -\dot\gamma_{n-1}/\gamma_{n-1}-z \cr  }\right],
$$
The diagonal elements of (24) yield $2\dot a_n = a_n(b_n-b_{n-1})$, and
the off-diagonal elements:
$2\dot\gamma_n/\gamma_n +b_n=0$,
$\dot a_n b_n+a_n\dot b_n = a_n(\dot\gamma_n/\gamma_n+\dot\gamma_{n-1}/
    \gamma_{n-1})b_n -2a_n^3 -a_nt$ and
$\dot a_n b_{n-1}+a_n\dot b_{n-1} = -a_n(\dot\gamma_n/\gamma_n+
\dot\gamma_{n-1}/\gamma_{n-1})b_{n-1} +2a_n^3 +a_nt$. Using (30), all
these equations are compatible with the differential system
$$\left\{ \matrix{ \dfrac{\dot a_n}{a_n} &=&
                                       b_n+\dfrac{n}{2a_n^2}, \cr
                           \dot b_n &=& -b_n^2  -2a_n^2 -t, \cr }
  \right. $$
which is the differential system equivalent to the second Painlev\'e
equation for $(-b_n,4a_n^2)$ ([Chu2], [Ge] p.339).
The connection with Painlev\'e transcendents can lead to advances in
the solution of the problem posed by Maroni in [Mar1]: when do we have
$a_1, a_2,\ldots \ne 0$ in (30)? The problem is now to localize the
zeros of solutions of special Painlev\'e equations.

{\bf 7. Example~3. $\exp(-x^4/4-tx^2)$ on {\blackb R}.}

This is the simplest nontrivial Freud's weight, and the corresponding
orthogonal polynomials have been much worked
 ([BoN] [Fr2] [LeQ] [Lub] [Mag2] [Mag3] [NeV] [GFOPCF] [NeV2] [Sho] ).
As for example~2,
we expand (13) and (15) with $W(z)=1$ and $2V(z)=-z^3-2tz$:
$$\Theta_n(z)=-z^2-2t-a_n^2-a_{n+1}^2,\qquad\qquad
  \Omega_n(z) = -z^3/2-(a_n^2+t)z. $$
A relation between the $a_n$'s is found by expanding (20), equating
the $z^2$ terms gives
$$ a_n^2(a_{n-1}^2+a_n^2+a_{n+1}^2)+2ta_n^2 =n ,\
     n=1,2,\ldots\  (a_0=0) \eqno(31) $$
a relation which seems to have been found by Shohat ([Sho]\  ),
rediscovered by Freud [Fr2] and Bessis [Bes]. Remark that we have a degree
of freedom on $a_1$: this is because the weight can have the real axis
{\it and\/} the pure imaginary axis in its support, with
$w(x)=\lambda\exp(-x^4/4-tx^2)$ on the pure imaginary axis, and the preceding
results hold for any $\lambda$, so $a_1$ is some function (which can be
computed from first moments) of $\lambda$. However, if it is
requested that all the $a_n$'s are positive, the solution is unique
and can be computed efficiently ([LeQ], see also [Nev2] p.470).
Now, (17) and (23) are computed:
$$A= \left[\matrix{ -a_n^2 z   &  a_n(z^2+2t+a_n^2+a_{n+1}^2)     \cr
         -a_n(z^2+2t+a_{n-1}^2+a_n^2) & z^3+(a_n^2+2t)z \cr  }\right],$$
$$
  H= \left[\matrix{\dot\gamma_n/\gamma_n-a_n^2  &     a_n z      \cr
               -a_n z &  -\dot\gamma_{n-1}/\gamma_{n-1}+z^2+a_n^2 \cr  }\right],
$$
The equations from (24) amount to be equivalent to
$$ \dfrac{\dot \gamma_n}{\gamma_n} = \dfrac{a_n^2+a_{n+1}^2}2 ,
  \eqno(32) $$
which, with $a_n\gamma_n=\gamma_{n-1}$, gives
$$ \dfrac{\dot a_n}{a_n} = \dfrac{a_{n-1}^2-a_{n+1}^2}2. \eqno(33)$$

{\it
Actually, $(32)$ (and $(33)$) can be recovered by quite
 elementary means:\/}
let $\{p_n(x;t)\}$ be the polynomials orthonormal with respect to an even
measure of the form
 $d\sigma(x;t)=\exp(-tx^2)d\sigma(x;0)$ on some support $S$,
we have then for the monic orthogonal polynomials $p_n/\gamma_n$:
$$\displaylines{
  \dfrac{\partial}{\partial t}\,\dfrac1{\gamma_n^2}=
  \dfrac{\partial}{\partial t}\,\int_S \left(
     \dfrac{p_n(x;t)}{\gamma_n}\right)^2 \exp(-tx^2)\,d\sigma(x;0) =
 \hfill\cr \hfill
= -\int_S x^2\left(\dfrac{p_n(x;t)}{\gamma_n}\right)^2
 \exp(-tx^2)\,d\sigma(x;0) = -\dfrac{a_n^2+a_{n+1}^2}{\gamma_n^2},\cr}$$
using $x^2p_n = a_na_{n-1}p_{n-2}+(a_n^2+a_{n+1}^2)p_n+a_{n+1}a_{n+2}
p_{n+2}$ from (2) when $b_n=0$, and that the derivative in $t$ of a
monic polynomial must be of degree $<n$. {\it Conversely\/}, it has
been shown that (33) implies that the $a_n$'s are the coefficients
of the recurrence of orthogonal polynomials with respect to a
measure of the form $\exp(-tx^2)d\sigma(x;0)$ where $d\sigma(x;0)$
does not depend on $t$ [KvM] [Mo] (see also  [Fra1], [Fra2]), [Y].

It is even probably possible to recover the information given by
(24) for all the semi-classical orthogonal polynomials by more
elementary means, but the connection with monodromy theory,
interesting on its own right, has more advantages: for
instance, it is known that the differential equations
produced by (24) have the Painlev\'e property (foreword of [Pain],
see [Mal] for a modern proof),
 i.e., movable singular points can only be poles (see [Cha], [In]
chap.\ 14). No wonder that the classical
Painlev\'e transcendants appear in these examples.

\noindent
Now, we get an equation for the single $a_n$ using (31): let
$u_n=a_n^2$, from $\dot u_n = u_n(u_{n-1}-u_{n+1})$,
$$\eqalign{
  \ddot u_n &= \dot u_n(u_{n-1}-u_{n+1})+u_n(\dot u_{n-1}-\dot u_{n+1})
         \cr
  &= u_n(u_{n-1}-u_{n+1})^2 +u_n[u_{n-1}(u_{n-2}-u_n) -
                                 u_{n+1}(u_n-u_{n+2})] \cr
  &= u_n(u_{n-1}-u_{n+1})^2 +u_n[n-1 -2tu_{n-1} -u_{n-1}(u_{n-1}+2u_n)+\cr
  &\hfill                       +n+1 -2tu_{n+1} -u_{n+1}(u_{n+1}+2u_n)]\cr
 &= u_n [ 2n -2(u_n+t)(u_{n-1}+u_{n+1}) -2u_{n-1}u_{n+1} ]   \cr
 &= u_n [ 2n -2(u_n+t)(u_{n-1}+u_{n+1}) -(u_{n-1}+u_{n+1})^2/2
                                        +(u_{n-1}-u_{n+1})^2/2 ]\cr
 &= u_n [2n+2(u_n+t)^2 -(u_{n-1}+2u_n+u_{n+1}+2t)^2/2]
        +(\dot u_n)^2/(2u_n) \cr
 &= u_n[ 2n+2(u_n+t)^2 -( n/{u_n}+u_n)^2/2] +(\dot u_n)^2/(2u_n) \cr
 &= \dfrac{u_n}2\left[ 4(u_n+t)^2 -\left(\dfrac n{u_n}-u_n\right)^2\right]
    +\dfrac{(\dot u_n)^2}{2u_n} \cr
 &= \dfrac{(\dot u_n)^2}{2u_n}  +\dfrac1{2u_n}
    \left(3u_n^2 +2tu_n -n\right)\left(u_n^2 +2tu_n +n\right),\cr
}$$
which is a special case of the $4^{\rm th}$ {\it Painlev\'e equation\/}
$$ \ddot y=\dfrac{\dot y^2}{2y}+\dfrac{3y^3}2+4ty^2 +2(t^2-\alpha)y
  +\dfrac{\beta}y \eqno(34)$$
with $\alpha=-n/2$ and $\beta=-n^2/2$
[Bu] [Fok1] [Fok2] [Ge] [Ok].

For $a_n=\sqrt{u_n}$, one has a form without first derivative:
$$ 4a_n^3 \ddot a_n = (3a_n^4 +2ta_n^2 -n)(a_n^4 +2ta_n^2+n).\eqno(35)$$

Let ${\gothi a}_n=n^{-1/4}a_n$ and ${\gothi t}=n^{-1/2}t$, then we
have another form
$$ 4{\gothi a}_n^3 \dfrac{d^2}{d{\gothi t}^2}{\gothi a}_n =
n^2 (3{\gothi a}_n^4 +2{\gothi{ta}}_n^2 -1)
    ({\gothi a}_n^4 +2{\gothi{ta}}_n^2+1).\eqno(36)$$

What can be the use of these equations? To explore these things, let us
first look at the graph of some ${\gothi a}_n$'s computed with
the Lew \& Quarles method [LeQ]:

\centerline{ \qquad ${\gothi a}_n=n^{-1/4}a_n$}

\vskip 5cm

\rightline{\qquad\qquad${\gothi t}=n^{-1/2}t$}

\vskip 5mm


\noindent
${\gothi a}_1,\ldots,{\gothi a}_{10}$ (only
           ${\gothi a}_1,\ldots,{\gothi a}_4$ are marked) tend to
be close to the zeros of the right-hand side of (36) (thick line).
In particular,
  ${\gothi a}_n(t)\sim 1/\sqrt{2{\gothi t}}$ when ${\gothi t}\to +\infty$:
$a_1^2=\mu_2/\mu_0 = \int_{-\infty}^\infty x^2 w(x)dx/
                             \int_{-\infty}^\infty  w(x)dx$, where
$w(x)=\exp(-x^4/4-tx^2)$. From [Erd] p.119,
$\mu_0 = \sqrt{\pi\sqrt2}\exp(t^2/2)D_{-1/2}(t\sqrt2)$ (parabolic
cylinder function). When $t\to +\infty$, $\mu_0\sim \sqrt{\pi/t}$
([Erd] p.122), so $a_1^2=\mu_2/\mu_0=-\dot\mu_0/\mu_0\sim 1/(2t)$.
From (31), if $a_1,a_2,\ldots, a_{n-1}$ are $O(t^{-1/2})$,
$a_n\sim \sqrt{n/(2t)}$ when $t\to\infty$.
When $t\to -\infty$, $\mu_0\sim$ constant $t^{-1/2}\exp(t^2)$
([Erd] p.123), so $a_1\sim\sqrt{-2t}$. The figure suggests that
${\gothi a}_n({\gothi t})\sim\sqrt{-2{\gothi t}}$ when $t\to -\infty$ and $n$
is odd, while
${\gothi a}_n({\gothi t})\sim\sqrt{-1/(2{\gothi t})}$ when $n$ is even.

One of the most interesting uses of expressions of recurrence coefficients
where $n$ is not bound to be an integer is to define
{\bf general associated orthogonal polynomials\/}, i.e., polynomials
defined by $a_{n+\nu+1}p_{n+1}^{(\nu)}(z)=
          (z-b_{n+\nu})p_{n}^{(\nu)}(z)-
            a_{n+\nu}p_{n-1}^{(\nu)}(z)$, and degree $p_n^{(\nu)}=n$
(as in [AW], [ILVW]).
So, let us define $a_{\nu}$ as some solution of (35) with $n$
replaced by $\nu$:
$$ 4a_{\nu}^3 \ddot a_{\nu} = (3a_{\nu}^4 +2ta_{\nu}^2 -\nu)
                              (a_{\nu}^4 +2ta_{\nu}^2+\nu),\eqno(37)$$
where $\nu$ is a given complex number. Then,
$y=\left[\dfrac{\nu}{2a_{\nu}^2}-\dfrac{a_\nu^2}2-t\mp\dfrac{\dot a_\mu}
       {a_\mu}\right]^{1/2}$ satisfies the same equation (37), but with
$\nu$ replaced by $\nu\pm1$ (Schlesinger transformation, [Fok1]
\S 3.3). Indeed, derivating
$y^2+t+(a_\nu^2-\nu/a_\nu^2)/2=\mp\dot a_\nu/a_\nu$ yields
$2y\dot y=\mp(\nu\pm 1-2a_\nu^2y^2-y^4-2ty^2)$ and a new derivation
establishes the property. So, the definition makes sense and
(31) still holds with $\nu$. There are still two degrees of freedom
in (37), but they are removed when suitable boundary conditions
are fixed ([DeC1] [DeC2]). Here, we just have to impose
$a_\nu =O(t^{-1/2)}$ when $t\to +\infty$ ([Yos], quoting Malmquist;
the point being that $a_\nu(t)$ must have an asymptotic series when
$t\to +\infty$ for fixed $\nu$, the relation with the dual situation,
i.e., $t$ fixed and $\nu\to +\infty$ is striking, see Section~4 of
[Wi]). In summary:

{\it For any real or complex $\nu$, the associated Freud orthogonal
  polynomials $p_n^{\nu}$ (which are related to the weight
  $\exp(-x^4/4-tx^2$ on {\blackb R} when $\nu=0$) have recurrence
 coefficient $a_{\nu+1}(t), a_{\nu+2}(t),\ldots$, where $a_\mu(t)$
 is completely defined as the solution of
$$ 4a_{\mu}^3 \ddot a_{\mu} = (3a_{\mu}^4 +2ta_{\mu}^2 -\mu)
                              (a_{\mu}^4 +2ta_{\mu}^2+\mu),$$
 which remains $O(t^{-1/2})$ when $t\to+\infty$.}

For the associated polynomials themselves, we can now construct
$\Theta_{\nu+n}$ and $\Omega_{\nu+n}$, therefore a differential
equation (18) with index $\nu+n$. Let $\varphi_{\nu+n}$ and
$\psi_{\nu+n}$ be two independent solutions of this differential
equation (in $z$).
 Following Hahn ([Ha1] eq. (17)),
$p_n^{(\nu)}=(\varphi_{\nu+n}\psi_{\nu-1}-\psi_{\nu+n}\varphi_{\nu-1})/
             (\varphi_\nu\psi_{\nu-1}-\psi_\nu\varphi_{\nu-1})$.
It can then be shown that
$f_\nu = \lim_{n\to\infty}p_{n-1}^{(\nu+1)}/p_n^{(\nu)}$
satisfies a {\it Riccati equation\/} (Laguerre-Hahn class [Mag1]).

{\bf 8. Example~4. $(x-t)^\rho \exp(-x^2)$ on $[t,\infty)$.}

The corresponding orthogonal polynomials are called (when $t=\rho=0$)
the Maxwell polynomials in [BeR], where other references can be
found ($\rho=1$: speed polynomials in [ClS]).
This case is closely related to the preceding one: put $x=t+u^2/2$ in
$\int_t^\infty p_n(x)p_m(x)(x-t)^\rho \exp(-x^2)\,dx=\delta_{m,n}$ to
find that $\tilde p_{2n}(u)=2^{-(\rho+1)/2}\exp(-t^2/2)p_n(t+u^2/2)$
is the orthonormal polynomial of degree $2n$ with respect to the
weight $\tilde w(u)=|u|^{2\rho+1}\exp(-u^4/4-tu^2)$ on {\blackb R}.
So, we have $a_n=\tilde a_{2n}\tilde a_{2n-1}/2$ and
$b_n = t+(\tilde a_{2n}^2+\tilde a_{2n+1}^2)/2$ ([Chi]   , etc.).

For the $\tilde a_n$'s, we still have $\dot{\tilde a}_n=
\tilde a_n(\tilde a_{n-1}^2-\tilde a_{n+1}^2)/2$ as before, but a
slightly different recurrence relation
$\tilde a_n^2(\tilde a_{n-1}^2+\tilde a_n^2+\tilde a_{n+1}^2+2t)=
n+(2\rho+1)\,$odd($n$), where odd($n)=(1-(-1)^n)/2$ [Fr2] [Mag2].
Working this yields now ($u_n = \tilde a_n^2)$.
$$\ddot u_n=\dfrac{\dot u_n^2}{2u_n} +\dfrac{3u_n^3}2 +4tu_n^2
 +2\left(t^2+\dfrac n2+(2\rho+1)\dfrac{1+3(-1)^n}4\right)u_n
 -\dfrac{(n+(2\rho+1){\rm odd}(n))^2}{2u_n},$$
i.e., the Painlev\'e $4^{\rm th}$ equation (34) with
$\alpha=-n/2-(2\rho+1)(1+3(-1)^n)/4$ and
$\beta=-(n+(2\rho+1){\rm odd}(n))^2/2$.

Many almost-classical orthogonal polynomials (see [Chin], [ClS]
and references in
[BeR] and [Gau]) could still be worked, and the simplest of them will
likely be related to other Painlev\'e transcendents (perhaps not the
{\it first\/} one\dots although [Fok2] finds first Painlev\'e
transcendents as solutions of a limit case of (31))
 At least a new case is briefly presented now:

{\bf 9. Example~5. Beyond Painlev\'e: $\exp(-x^6-tx^2)$ on {\blackb R}.}

With $u_n=a_n^2$, $\dot u_n=u_n(u_{n-1}-u_{n+1})$ still holds, but
the recurrence relation is somewhat more complicated than before
[Fr2] [Mag2] [Mag3]:
$u_n(u_{n-2}u_{n-1}+u_{n-1}^2+2u_{n-1}u_n +u_n^2 +2u_nu_{n+1}
 +u_{n-1}u_{n+1}+u_{n+1}^2+u_{n+1}u_{n+2}+2t)=n$, for
$n=1,2,\ldots$ As a first step, one has a differential system for
$u_{n-1},\ldots,u_{n+2}$ by eliminating $u_{n-2}$ and $u_{n+3}$
from the recurrence relation:
$$\left\{
  \eqalign{
 \dot u_{n-1} &= u_{n-1}u_{n-2}-u_{n-1}u_n= \cr
    &=\dfrac n{u_n}-2t-u_{n-1}^2-3u_{n-1}u_n-u_n^2-2u_nu_{n+1}
           -u_{n-1}u_{n+1}-u_{n+1}^2-u_{n+1}u_{n+2},\cr
 \dot u_n&= u_n(u_{n-1}-u_{n+1}),\cr
 \dot u_{n+1}&= u_{n+1}(u_n-u_{n+2}),\cr
 \dot u_{n+2} &= u_{n+1}u_{n+2}-u_{n+2}u_{n+3}= \cr
&=-\dfrac{n+1}{u_{n+1}}+2t+u_{n+2}^2+u_nu_{n+2}+3u_{n+1}u_{n+2}+
                                               u_{n+1}^2+2u_nu_{n+1}
           +u_n^2+u_{n-1}u_n,\cr
           } \right.
$$
which can still be transformed\dots This case in considered in [Fok2].

\vfill

{\bf Acknowledgements.}

Many thanks to R.Askey, R.Caboz,
D.\& G.Chudnovsky, C. De Coster,
L.Haine, J.Meinguet, A.Ronveaux,
P.\ van Moerbeke, M.Willem.

\vfill\eject

{\bf References.}

    \parskip=12pt   \baselineskip=10pt
\frenchspacing      \parindent=40pt

\item{[Ak]} N.I.~AKHIEZER, {\sl Elements of the Theory of Elliptic
          Functions\/}, translated from the $2^{\rm nd}$ Russian
          edition (Nauka, Moscow, 1970), {\sl Transl.\ Math.\ Monographs\/}
          {\bf 79\/}, A.M.S., Providence, 1990.

\item{[Al]} I.\ ALVAREZ ROCHA, F.MARCELLAN.
             On semiclassical linear functionals:
        integral representations. This volume.

\item{[Apt]}  A.I.\ APTEKAREV,  Asymptotic properties of
polynomials orthogonal on a system of contours and periodic motions of
Toda lattices, {\sl Mat.\ Sb.\ } {\bf 125\/} (1984) 231-258, =
 {\sl Math. USSR Sbornik\/} {\bf 53\/} (1986) 233-260.

\item{[AW]} R.ASKEY, J.WIMP, Associated Laguerre and Hermite polynomials,
     {\sl Proc.\ Royal Soc.\ Edinburgh\/} {\bf 96A\/} (1984) 15-37.

\item{[AtE]} F.V.ATKINSON, W.N.EVERITT, Orthogonal polynomials which
       satisfy second order differential equations, pp.~173-181
       {\it in\/} {\sl E.B.Christoffel\/} (P.L.BUTZER and F.FEH\'ER,
       editors), Birkh\"auser, Basel, 1981.

\item{[Bel]} S.\ BELMEHDI, On semi-classical linear functionals of
        class $s=1$. Classification and integral representations.
       {\sl Indag.\ Mathem.\ \/} N.S.\ {\bf 3\/} (3) (1992),
       253-275.

\item{[BeR]}  S. BELMEHDI, A. RONVEAUX, Laguerre-Freud's equations
     for  the recurrence coefficients of semi-classical
     orthogonal polynomials, to appear in {\sl J.\ Approx.\ Theory\/}.

\item{[BeR2]}  S. BELMEHDI, A. RONVEAUX, On the coefficients of the
     three-term recurrence relation satisfied by some orthogonal
     polynomials. {\sl Innovative Methods in Numerical Analysis\/},
     Bressanone, Sept.\ 7-$11^{\rm th}$, 1992.

\item{[Bes]} D.~BESSIS, A new method in the combinatorics of the
    topological expansion, {\sl Comm.\ Math.\ Phys.\/} {\bf 69\/}
    (1979), 147-163.

\item{[BIZ]} D.\ BESSIS, C.\ ITZYKSON, J.B.\ ZUBER, Quantum field
     theory techniques in graphical enumeration, {\sl Adv.\ in
     Appl.\ Math.\/} {\bf 1\/} (1980), 109-157.

\item{[BoN]}  S. BONAN, P.NEVAI, Orthogonal polynomials and their
       derivatives,I, {\sl J. Approx. Theory\/} {\bf 40\/} (1984),
       134-147.

\item{[BLN]}  S.S. BONAN, D.S. LUBINSKY, P.NEVAI, Orthogonal polynomials
 and tand their
       derivatives,II, {\sl SIAM J. Math. An.\/} {\bf 18\/} (1987),
       1163-1176.

\item{[Bre]} C.~BREZINSKI, {\sl Pad\'e-type Approximation and General
          Orthogonal Polynomials \ ISNM\/} {\bf 50\/},
          Birkh\"auser-Verlag, Basel, 1980.

\item{[Brez]} C.~BREZINSKI, {\sl History of Continued Fractions and Pad\'e
           Approximants\/}, Springer-Verlag, Berlin, 1991.

\item{[dBvR]} M.G.~de~BRUIN, H.~van~ROSSUM, Formal Pad\'e approximation,
             {\sl Nieuw Arch. Wisk. $(3)$\/} {\bf 23\/} (1975), 115-130.

\item{[Bu]} F.J.~BUREAU, Les \'equations diff\'erentielles du second
      ordre \`a points critiques fixes, II. Les int\'egrales de
      l'\'equation A4 de Painlev\'e, {\sl Bull.\ Cl.\ Sci.\ Acad.\
      Roy.\ Belg.\/} {\bf 69\/} (1983) 397-433.

\item{[Cha]} R.\ CHALKLEY, New contributions to the related works of
      Paul Appell, Lazarus Fuchs, Georg Hamel, and Paul Painlev\'e on
     nonlinear differential equations whose solutions are free of
     movable branch points. {\sl J.\ Diff.\ Eq.\ \/} {\bf 68\/} (1987)
     72-117.

\item{[Chi]} T.S.~CHIHARA, ``An Introduction to Orthogonal Polynomials,''
               Gordon \& Breach, New York, 1978.

\item{[Chin]}  R.C.Y.~CHIN, A domain decomposition method for generating
            orthogonal polynomials for a Gaussian weight on a finite
            interval, {\sl J. Comp. Phys.\/} {\bf 99\/} (1992) 321-336.

\item{[Chua]} D.V.~CHUDNOVSKY, Riemann monodromy problem, isomonodromy
      deformation eqations and completely integrable systems, pp.385-447
      {\it in\/} {\sl Bifurcation Phenomena in Mathematical Physics and
      Related Topics, Proceedings  Carg\`ese, 1979\/} (C.BARDOS \&
      D.BESSIS, editors), NATO ASI series~C, vol. {\bf 54\/}, D.Reidel,
      Dordrecht, 1980.

\item{[Chub]} G.V.~CHUDNOVSKY, Pad\'e approximation and the
      Riemann monodromy problem, pp.449-510
      {\it in\/} {\sl Bifurcation Phenomena in Mathematical Physics and
      Related Topics, Proceedings  Carg\`ese, 1979\/} (C.BARDOS \&
      D.BESSIS, editors), NATO ASI series~C, vol. {\bf 54\/}, D.Reidel,
      Dordrecht, 1980.

\item{[Chu0]} D.V.CHUDNOVSKY, G.V.CHUDNOVSKY, Introduction to
                           {\sl The Riemann Problem, Complete
     Integrability and Arithmetic Applications\/} (D.~Chudnovsky and
     G.~Chudnovski, Eds.), pp.1-11, Springer-Verlag (Lecture
     Notes Math. {\bf 925\/}), Berlin, 1982.

\item{[Chu1]} D.V.CHUDNOVSKY, G.V.CHUDNOVSKY, Laws of composition of
   B\"acklund transformations and the universal form of completely
   integrable systems in dimensions two and three, {\sl Proc.\ Nat.\
    Acad.\ Sci.\ USA\/} {\bf 80\/} (1983) 1774-1777.

\item{[Chu2]} D.V.CHUDNOVSKY, G.V.CHUDNOVSKY, High precision computation
   of special function in different domains, talk given at Symbolic
    Mathematical Computation conference, Oberlech, July 1991.

\item{[ClS]} A.S.\ CLARKE, B.\ SHIZGAL, On the generation of orthogonal
      polynomials using asymptotic methods for recurrence coefficients.
      {\sl J.\ Comp.\ Phys.\ \/} {\bf 104\/} (1993) 140-149.

\item{[DeC1]} C.DE~COSTER, M.WILLEM, private communication,
    18 September 1992.

\item{[DeC2]} C.DE~COSTER, M.WILLEM, Density, spectral theory and
      homoclinics for singular Sturm-Liouville systems, to appear
     in {\sl J.\ Comp.\ Appl.\ Math.}

\item{[Dra]} A.~DRAUX, {\sl Polyn\^omes orthogonaux formels $-$
        Applications\/}, Lect.\ Notes Math.\ {\bf 974\/}, Springer-Verlag,
        Berlin, 1983.

\item{[Erd]} A.~ERD\'ELYI {\it et al.\/}, editors, {\sl Higher Transcendental
          Functions\/}, vol.~II, McGraw-Hill, New York, 1953.

\item{[Fed]} M.V.~FEDORYUK, Isomonodromy deformations of equations with
        irregular singularities, {\sl Mat.\ Sb.\/} {\bf 181\/} (1990) =
        {\sl Math.\ USSR Sb.\/} {\bf 71\/} (1992) 463-479.

\item{[Fok1]} A.S.\ FOKAS, U.\ MUGAN, M.J.\ ABLOWITZ, A method of
         linearization for Painlev\'e equations. Painlev\'e IV,V.
         {\sl Physica D\/} {\bf 30\/} (1988) 247-283.

\item{[Fok2]} A.S.\ FOKAS, A.R.\ ITS, A.V.\ KITAEV, Discrete
        Painlev\'e equations and their appearance in quantum gravity,
       {\sl Commun.\ Math.\ Phys.\ \/} {\bf 142\/} (1991) 313-344.

\item{[Fra1]} J.P.~FRANCOISE, Symplectic geometry and integrable
    $m$-body problems on the line, {\sl J.\ Math.\ Phys.\/}{\bf 29\/}
    (1988) 1150-1153.

\item{[Fra2]} J.P.~FRANCOISE, Syst\`emes int\'egrables \`a $m$ corps
      sur la droite, {\it in\/} {\sl Analyse Globale et Physique
      Math\'ematique\/}, Dec. 1989, preprint.

\item{[Fr1]}  G.~FREUD, {\sl Orthogonal Polynomials\/}, Akad\'emiai
            Kiad\'o/Pergamon Press, Budapest/Oxford, 1971.

\item{[Fr2]}  G.FREUD, On the coefficients in the recursion formul\ae\  of
     orthogonal polynomials, {\sl Proc.\ Royal Irish Acad. Sect. A\/}
     {\bf 76\/} (1976), 1-6.

\item{[RFu]} R.\ FUCHS, \"Uber lineare homogene Differentialgleichungen
       zweiter Ordnung mit drei im Endlichen gelegenen wesentlich
      singul\"are Stellen, {\sl Math.\ Ann.\/} {\bf 63\/} (1907)
      301-321.

\item{[GaN]} J.L.~GAMMEL, J.~NUTTALL, Note on generalized Jacobi
     polynomials, {\it in\/} ``The Riemann Problem, Complete
     Integrability and Arithmetic Applications'' (D.~Chudnovsky and
     G.~Chudnovski, Eds.), pp.258-270, Springer-Verlag (Lecture
     Notes Math. {\bf 925\/}), Berlin, 1982.

\item{[GaCL]} J.P.\ GASPARD, F.\ CYROT-LACKMANN, Density of states from
     moments. Application to the impurity band, {\sl J.\ Phys.\ C:
      Solid State Phys.\/} {\bf 6\/} (1973) 3077-3096.

\item{[Gau]} W.GAUTSCHI, Computational aspects of orthogonal
    polynomials, pp.181-216 {\it in\/} {\sl Orthogonal Polynomials:
    Theory and Practice\/} (P.NEVAI, editor) {\sl NATO ASI Series~C\/}
    {\bf 294\/}, Kluwer, Dordrecht, 1990.

\item{[Ge]} R.GERARD, La g\'eom\'etrie des transcendantes de
     P.Painlev\'e, pp.323-352 {\it in\/} {\sl Math\'ematique et
    Physique, S\'eminaire de l'Ecole Normale Sup\'erieure 1979-1982\/}
    (L.~BOUTET de MONVEL, A.DOUADY \& J.L.VERDIER, editors),
    {\sl Progress in Mathematics\/} {\bf 37\/}, Birkh\"auser,
    Boston, 1983.

\item{[GV1]} J.S.GERONIMO, W.VAN~ASSCHE, Orthogonal polynomials with
    asymptotically periodic recurrence coefficients, {\sl J.\ Approx.\
    Th.\/} {\bf 46\/} (1986) 251-283.

\item{[GV2]} J.S.GERONIMO, W.VAN~ASSCHE, Approximating the weight
    function for orthogonal polynomials on several intervals,
    {\sl J.\ Approx.\ Th.\/} {\bf 65\/} (1991), 341-371.

\item{[Gr]} C.C.~GROSJEAN, The measure induced by orthogonal polynomials
       satisfying a recursion formula with either constant or periodic
       coefficients. Part~I: Constant coefficients, {\sl Acad.\
       Analecta, Kon.\ Acad.\ Wet.\ Lett.\ Sch.\ Kunsten Belg.\/}
       {\bf 48\/}, Nr.~3, (1986), 39-60. Part~II: Pure or mixed
       periodic coefficients (general theory), {\it ibid.\/} {\bf 48\/}
       Nr.~5 (1986) 55-94.

\item{[GrM1]} D.J.\ GROSS, A.A.\ MIGDAL, A nonperturbative treatment
       of two-dimensional quantum gravity. Princeton preprint PUPT
       1159 (1989).

\item{[GrM2]} D.J.\ GROSS, A.A.\ MIGDAL, Nonperturbative
        two-dimensional quantum gravity. {\sl Phys.\ Rev.\ Letters\/}
       {\bf 64\/} (1990) 127-130.

\item{[Ha1]} W.HAHN, On differential equations for orthogonal
      polynomials, {\sl Funk. Ekvacioj\/}, {\bf 21\/} (1978) 1-9.

\item{[Ha2]} W.HAHN, \"Uber Orthogonalpolynome, die linearen
      funktionalgleichungen gen\"ugen, pp. 16-35 {\it in\/}
      {\sl Polyn\^omes Orthogonaux et Applications, Proceedings,
      Bar-le-Duc 1984\/}, (C.BREZINSKI \& al., editors),
      {\sl Lecture Notes Math.\/} {\bf 1171\/}, Springer, Berlin 1985.

\item{[HH]} L.\ HAINE, E.\ HOROZOV, Toda orbits of Laguerre
      polynomials and representations of the Virasoro algebra.
      Preprint Institut Math\'ematique Universit\'e Catholique de
     Louvain 217 (1992).

\item{[Hay]} R.\ HAYDOCK, The recursive solution of the Schr\"odinger
    equation, pp.215-294 {\it in\/} H.\ EHRENREICH {\it et al.\/},
    editors: {\sl Solid State Physics\/} {\bf 35\/} , Ac.\ Press,
    N.Y.\ , 1980.

\item{[HayN]} R.\ HAYDOCK, C.M.M.\ NEX, A general terminator for the
       recursion method, {\sl J.\ Phys.\ C: Solid State Phys.\ }
       {\bf 18\/} (1985) 2235-2248.

\item{[HvR1]}  E. HENDRIKSEN, H. van ROSSUM, A Pad\'e-type approach to
     non-classical orthogonal polynomials, {\sl J. Math. An. Appl.\/}
     {\bf 106\/} (1985) 237-248.

\item{[HvR2]}  E. HENDRIKSEN, H. van ROSSUM, Semi-classical orthogonal
     polynomials,
                                       pp. 354-361 {\it in\/}
      {\sl Polyn\^omes Orthogonaux et Applications, Proceedings,
      Bar-le-Duc 1984\/}, (C.BREZINSKI \& al., editors),
      {\sl Lecture Notes Math.\/} {\bf 1171\/}, Springer, Berlin 1985.

\item{[In]} E.L.~INCE, {\sl Ordinary Differential Equations\/},
      Longmans Green 1928 = Dover 1956.

\item{[I]} M.ISMAIL, On sieved orthogonal polynomials III: orthogonality
        on several intervals, {\sl Trans.\ Amer.\ Math.\ Soc.\/}
        {\bf 294\/} (1986), 89-111.

\item{[ILVW]} M.ISMAIL, J.LETESSIER, G.VALENT, J.WIMP, Some results
      on associated Wilson polynomials, pp.293-298 {\it in\/}
      {\sl Orthogonal Polynomials and their Applications\/}
      (C.BREZINSKI {\it et al.\/}, editors),
      {\sl IMACS Annals on Computing and Applied Mathematics\/}
      {\bf 9\/} (1991), Baltzer AG, Basel.

\item{[KvM]} M.KAC, P.~van~MOERBEKE, On an explicitly soluble system
     on nonlinear differential equations related to certain Toda
     lattices, {\sl Adv.\ Math.\/} {\bf 16\/} (1975) 160-169.

\item{[Lag]}  E. LAGUERRE, Sur la r\'eduction en fractions continues d'une
      fraction qui satisfait \`a une \'equation diff\'erentielle lin\'eaire
      du premier ordre dont les coefficients sont rationnels,
      {\sl J. Math. Pures Appl. (4)\/} {\bf 1\/} (1885), 135-165 =
      pp. 685-711 {\it in\/}
      {\sl Oeuvres\/}, Vol.II, Chelsea, New-York 1972.

\item{[LaG]} Ph.\ LAMBIN, J.P.\ GASPARD, Continued-fraction technique for
      tight-binding systems. A generalized-moments approach,
      {\sl Phys.\ Rev.\ B/} {\bf 26\/} (1982) 4356-4368.

\item{[LD]}  J.A.LAPPO-DANILEVSKY, {\sl M\'emoires sur la th\'eorie
        des syst\`emes des \'equa\-tions diff\'erentielles lin\'eaires\/},
        vol.~I,II,III bound as one volume, Chelsea Pub.\ Co.\ ,1953.

\item{[LW]} D.\ LEVI, P.\ WINTERNITZ, editors: {\sl Painlev\'e
       Transcendents. Their Asymptotics and Physical Applications.
       NATO ASI Series: Series~B: Physics\/} {\bf 278\/},
       Plenum Press, N.Y., 1992.

\item{[LeQ]} J.S.LEW, D.A.QUARLES, Nonnegative solutions of a nonlinear
     recurrence, {\sl J.\ Approx.\ Th.\/} {\bf 38\/} (1983), 357-379.

\item{[LiMu]} J.-M.\ LIU, G.\ M\"ULLER, Infinite-temperature dynamics of
      the equivalent-neighbor $XYZ$ model, {\sl Phys.\ Rev.\ A\/}
     {\bf 42\/} (1990) 5854-5864.

\item {[Lub]} D.S.~LUBINSKY, A survey of general orthogonal polynomials
     for weights on finite and infinite intervals, {\sl Acta Applicand\ae
     Mathematic\ae\/} {\bf 10\/} (1987) 237-296.

\item  {[Mag1]}  A.P.~MAGNUS, Riccati acceleration of Jacobi
      continued fractions and Laguerre-Hahn orthogonal polynomials,
     pp.~213-230 {\it in\/} {\sl Pad\'e Approximation and its
       Applications, Proceedings, Bad Honnef 1983\/}, Lecture Notes
     Math. {\bf 1071\/} (H.WERNER \& H.T.B\"UNGER, editors),
      Springer-Verlag, Berlin, 1984.

\item  {[Mag2]}  A.P.~MAGNUS, A proof of Freud's conjecture about
     orthogonal polynomials related to $|x|^\rho \exp(-x^{2m})$ for
     integer $m$, pp.~362-372 {\it in\/} {\sl Polyn\^omes Orthogonaux
     et Applications, Proceedings Bar-le-Duc 1984.\/} (C.BREZINSKI
     {\it et al.\/}, editors), {\sl Lecture Notes Math.\/} {\bf 1171\/},
     Springer-Verlag, Berlin 1985.

\item  {[Mag3]}  A.P.~MAGNUS, On Freud's equations for exponential
     weights, {\sl J.\ Approx.\ Th.\/} {\bf 46\/} (1986) 65-99.

\item  {[Mag4]}  A.P.~MAGNUS, Associated Askey-Wilson polynomials as
     Laguerre-Hahn orthogonal polynomials, pp.~261-278 {\it in\/}
    {\sl Orthogonal Polynomials and their Applications, Proceedings
     Segovia 1986.\/} (M.ALFARO {\it et al.\/}, editors), {\sl Lecture
     Notes Math.\/} {\bf 1329\/}, Springer-Verlag, Berlin 1988.

\item{[Mal]} B.MALGRANGE, Sur les d\'eformations isomonodromiques.
    I. Sigularit\'es r\'eguli\`eres,
                   pp.401-426 {\it in\/} {\sl Math\'ematique et
    Physique, S\'eminaire de l'Ecole Normale Sup\'erieure 1979-1982\/}
    (L.~BOUTET de MONVEL, A.DOUADY \& J.L.VERDIER, editors),
    {\sl Progress in Mathematics\/} {\bf 37\/}, Birkh\"auser,
    Boston, 1983. II. Singularit\'es irr\'eguli\`eres, {\it ibid.\/},
    427-438.

\item {[Mar]}  P. MARONI, Le calcul des formes lin\'eaires et les
      polyn\^omes orthogonaux semi-classiques,
      pp. 279-290 {\it in\/} {\sl Orthogonal Polynomials and their
      Applications. Proceedings, Segovia 1986\/} (M. ALFARO
      {\it et al.\/}, editors), {\sl Lecture Notes Math.\/}
      {\bf 1329\/}, Springer, Berlin 1988.

\item {[Mar1]}  P. MARONI, Un exemple d'une suite orthogonale
        semi-classique de classe un. {\sl Polinomios ortogonales y
      applicaciones; Actas VI Simposium, Gijon\/} (1989), 234-241.

\item {[Mar2]}  P. MARONI, Une th\'eorie alg\'ebrique des
      polyn\^omes orthogonaux. Application aux
      polyn\^omes orthogonaux semi-classiques. pp. 95-130 {\it in\/}
      {\sl Orthogonal Polynomials and their Applications\/}
      (C.BREZINSKI {\it et al.\/}, editors),
      {\sl IMACS Annals on Computing and Applied Mathematics\/}
      {\bf 9\/} (1991), Baltzer AG, Basel.

\item{[Mo]} J.MOSER, Three integrable Hamiltonian systems connected
       with isospectral deformations, {\sl Adv.\ Math.\/} {\bf 16\/}
      (1975) 197-220.

\item{[Nev]} P.NEVAI, Two of my favorite ways of obtaining asymptotics
    for orthogonal polynomials, pp.417-436 {\it in\/}
    {\sl Anniversary Volume on Approximation Theory and
      Functional Analysis\/}, (P.L.BUTZER, R.L.STENS and B.Sz.-NAGY,
      editors), {\bf ISNM 65\/}, Birkh\"auser Verlag, Basel, 1984.

\item {[GFOPCF]}  P. NEVAI, G\'eza Freud, orthogonal polynomials and
     Christoffel functions. A case study, {\sl J. Approx. Theory\/}
     {\bf 48\/} (1986), 3-167.

\item{[Nev2]} P.NEVAI, Research problems in orthogonal polynomials,
    pp. 449-489 {\it in\/} {\sl Approximation Theory VI\/}, vol.
   {\bf 2\/} (C.K.CHUI, L.L.SCHUMAKER \& J.D.WARD, editors),
    Academic Press, 1989.

\item{[N]} J.NUTTALL, Asymptotics of diagonal Hermite-Pad\'e polynomials,
 {\sl J.Approx.\ Th.\ } {\bf 42\ } (1984) 299-386.

\item{[Ok]} K.\ OKAMOTO, Studies on the Painlev\'e equations III.
        Second and fourth Painlev\'e equations,
       ${\rm P}_{\rm II}$ and
       ${\rm P}_{\rm IV}$, {\sl Math.\ Ann.\/} {\bf 275\/} (1986)
       221-255.

\item{[OW]} E.P.\ O'REILLY, D.\ WEAIRE, On the asymptotic form  of the
      recursion basis vectors for periodic Hamiltonians, {\sl J.\ Phys.\ A:
      Math.\ Gen.\/} {\bf 17\/} (1984) 2389-2397.

\item{[Pain]} P.\ PAINLEV\'E, {\sl \OE uvres de Paul Painlev\'e\/},
     vol.\ {\bf 3\/}, C.N.R.S.\ , Paris, 1975.

\item {[Peh1]} F.PEHERSTORFER, On Bernstein-Szeg\H o orthogonal
      polynomials on several intervals. II. Orthogonal polynomials
      with periodic recurrence coefficients, {\sl J.\ Approx.\ Th.\/}
      {\bf 64\/} (1991) 123-161.

\item {[Peh]} F.PEHERSTORFER, On orthogonal polynomials on several
       intervals, {\sl VII Simposium sobre polinomios ortogonales
       y applicaciones\/}, Granada, Espa\~na, 23-27 Sept.\ 1991.

\item {[Per]} O.PERRON, {\sl Die Lehre von den Kettenbr\"uchen\/},
        $2^{\rm nd}$ edition, Teubner, Leipzig, 1929 = Chelsea,

\item {[Sho]}  J.A. SHOHAT, A differential equation for orthogonal
     polynomials, {\sl Duke Math. J.\/} {\bf 5\/} (1939),401-417.

\item{[StT]} H.STAHL, V.TOTIK, {\sl General Orthogonal Polynomials\/},
         ({\sl Encyc.\ Math.\ Appl.\/} {\bf 43\/}), Cambridge U.P.,
         Cambridge, 1992.

\item{[VA]} W.\ VAN~ASSCHE, {\sl Asymptotics for Orthogonal
      Polynomials. Springer Lecture Notes Math.\/} {\bf 1265\/},
     Springer-Verlag, Berlin 1987.

\item{[Wi]} J.WIMP, Current trends in asymptotics: some problems and
      some solutions, {\sl J.\ Comp.\ Appl.\ Math.\/} {\bf 35\/} (1991)
      53-79.

\item{[Y]} S.YAMAZAKI, The semi-infinite system of nonlinear
     differential equations $\dot A_k = 2A_k(A_{k+1}-A_{k-1})$;
    methods of integration and asymptotic time behaviours,
    {\sl Nonlinearity\/} {\bf 3\/} (1990) 653-676.

\item{[Yos]} S.YOSHIDA, 2-parameter family of solutions for Painlev\'e
      (I)$\sim$(V) at an irregular singular point,
                   {\sl Funk. Ekvacioj\/}, {\bf 28\/} (1985) 233-248.

\item{[Zu]} J.B.\ ZUBER, L'invariance conforme et la physique \`a deux
       dimensions. {\sl La Recherche\/} {\bf 24\/} (1993) 142-151.

\bye